\documentclass[11pt,a4paper]{amsart}
\usepackage[left=2cm, right=2cm, top=2.8cm, bottom=2.8cm]{geometry}

\allowdisplaybreaks
\usepackage{booktabs}

\usepackage{tikz}
\usepackage{pgfplots}
\pgfplotsset{compat=1.18}
\usetikzlibrary{arrows.meta,calc}

\usepackage{algorithm}
\floatname{algorithm}{Algorithm}
\usepackage{subcaption}
\usepackage[colorlinks=true,citecolor=red,linkcolor=blue,urlcolor=cyan,filecolor=magenta]{hyperref}
\usepackage{lmodern}            % Latin Modern font
\usepackage[T1]{fontenc}        % Font encoding
\usepackage[utf8]{inputenc}     % Input encoding
\usepackage{indentfirst}        % Indent the first paragraph of each section
\usepackage{color,xcolor}       % Colors
\usepackage{graphicx}           % Graphics
\usepackage{microtype}          % Improved justification
\usepackage{amsbsy}             % Bold mathematical symbols
\usepackage{amssymb}            % AMS symbols
\usepackage{amsthm}             % Theorem environments
\usepackage{mathrsfs}           % \mathscr
\usepackage{lineno}
\usepackage{tikz-cd}
\usepackage{comment}
\usetikzlibrary{trees}

\DeclareMathOperator{\dist}{dist}
\DeclareMathOperator{\conv}{co}
\DeclareMathOperator{\Prox}{Prox}
\DeclareMathOperator{\intr}{int}

\newtheorem{theorem}{Theorem}[section]
\newtheorem{corollary}[theorem]{Corollary}
\newtheorem{lemma}[theorem]{Lemma}
\newtheorem{proposition}[theorem]{Proposition}
\newtheorem{remark}[theorem]{Remark}

\newtheorem{example}[theorem]{Example}

\title[Strong quasiconvexity with the exact modulus]{First-order characterizations of strong quasiconvexity with the exact modulus}

\author[J. Matos Ascona]{John Matos Ascona}
\address{Universidade Estadual de Mato Grosso do Sul, Rodovia Dourados-Itahum (Rodovia MS/162), km 12, Caixa Postal 351, Cidade Universit\'aria de Dourados, 79804-970, Dourados, MS, Brazil}
\email{john.ascona@uems.br}

\author[J. C. Huaccha Neyra]{Jackeline C. Huaccha Neyra}
\address{Instituto de Ci\^encia e Tecnologia, Universidade Federal de S\~ao Paulo,
Av. Cesare Mansueto Giulio Lattes, 1201, 12247-014 S\~ao Jos\'e dos Campos, SP, Brazil}
\email{jackeline.neyra@unifesp.br}

\author[L. Santa Cruz Calder\'on]{Lizet Santa Cruz Calder\'on}
\address{Universidade Federal de Mato Grosso do Sul, Rua Oscar Trindade de Barros, 740, Bairro da Serraria, 79200-000, Aquidauana, MS, Brazil}
\email{lizet.calderon@ufms.br}
\thanks{Corresponding author: Lizet Santa Cruz Calder\'on (lizet.calderon@ufms.br).}

\author[H. L\'opez-L\'azaro]{Heraclio L\'opez-L\'azaro}
\address{Departamento de Matem\'atica, Universidade Federal de Santa Catarina,
Campus Universit\'ario Reitor Jo\~ao David Ferreira Lima, Trindade, 88040-900 Florian\'opolis, SC, Brazil}
\email{heraclio.lopez.lazaro@ufsc.br}

\subjclass[2020]{90C26, 49J52, 26B25, 26A24, 65K10, 34A60}
\keywords{Strongly quasiconvex functions, Dini derivatives, Clarke subdifferential, exact modulus, subgradient method}

\begin{document}
\allowdisplaybreaks % allow page breaks inside displayed (in)equalities

\begin{abstract}
In this paper, we study first-order characterizations of strong quasiconvexity for functions defined on convex subsets of real normed spaces. For functions that are continuous along segments, an inequality involving the lower Dini derivative is a known necessary condition for strong quasiconvexity with a given modulus, whereas the established converse guarantees only half of that modulus. By means of a one-dimensional lemma, we prove that this inequality, and even its weaker version with the upper Dini derivative, ensures strong quasiconvexity with the same modulus, without requiring the domain to be open; this answers a question raised in the recent literature. The main result leads to a characterization with the same modulus in terms of the Clarke subdifferential for locally Lipschitz functions on Banach spaces, which answers another open question, and to a formula for the exact modulus of radial functions defined on Hilbert spaces. We also present some consequences in optimization, such as an optimality condition without Clarke regularity, an error bound with optimal constant and, for the projected subgradient method, the same convergence rate as in the strongly convex case, together with the exponential convergence of the associated flow. Finally, we show that the converse of a relation between proximal points and strong subdifferentials fails in general, but holds for a truncation of the function.
\end{abstract}
\maketitle

\section{Introduction}\label{Sec. Introd}
One of the most natural extensions of strong convexity to nonconvex optimization is the class of strongly quasiconvex functions, introduced by Polyak \cite{Polyak1966} as a particular case of uniformly quasiconvex functions, see \cite[Section 1]{grad2025strongly}. Let $X$ be a real normed space and let $C\subseteq X$ be a convex set. In the form in which it is currently used, cf. \cite{grad2025strongly}, a function $h:C\to \mathbb{R}$ is strongly quasiconvex with modulus $\gamma>0$ if
\begin{equation}\label{eq. SQ}
h (x+t(y-x))\leq\max\{h(x),h(y)\}-\frac{\gamma}{2}t(1-t)\|y-x\|^2\quad \text{for all }x,y\in C,~ t\in[0,1].
\end{equation}
If $C$ has more than one point and $h$ satisfies \eqref{eq. SQ} with some modulus, then the set of admissible moduli is an interval of the form $(0,\gamma^*]$, with $\gamma^*<\infty$, and $\gamma^*$ is called the exact modulus of $h$ on $C$ (see also the notion of exact modulus of quasiconvexity in \cite[Definition 1]{vladimirov1978uniformly}).

Every strongly convex function with modulus $\gamma$ satisfies \eqref{eq. SQ}, but the class is broader and contains functions that are neither convex nor locally Lipschitz (see Example \ref{Ex. Coro radial}). When $X=\mathbb{R}^n$, these functions are $2$-supercoercive and, if $C$ is closed and $h$ is lower semicontinuous, then $h$ has exactly one minimizer on $C$ \cite[Theorem 1 and Corollary 3]{lara2022strongly}. For this reason, they constitute a natural setting for proximal, gradient and subgradient methods beyond convexity \cite{choque2024subgradient,lara2022strongly,lara2025characterizations}.

The first-order characterization of this class goes back to Vladimirov, Nesterov and Chekanov \cite[Theorems 1 and 6]{vladimirov1978uniformly} (see also \cite[Theorem 1]{jovanovic1993strongly}), according to whom, for a function $h$ of class $C^1$, strong quasiconvexity with modulus $\gamma$ implies
\begin{equation}\label{eq. grad}
h(x)\leq h(y) \implies  \langle\nabla h(y),y-x\rangle\geq\frac{\gamma}{2}\|y-x\|^2\quad \text{for all }x,y\in C,
\end{equation}
and, conversely, \eqref{eq. grad} implies strong quasiconvexity with modulus $\frac{\gamma}{2}$. The equivalence with the same modulus $\gamma$ in both directions was established by Bao and Nam \cite[Theorem 3.2]{bao2025first} for G\^ateaux differentiable functions on open sets.

More recently, in the nonsmooth setting, Hadjisavvas and Lara \cite[Theorem 6(b)]{hadjisavvas2026characterizations} showed that, for $h$ continuous along segments of an open convex set, strong quasiconvexity with modulus $\gamma$ implies
\begin{equation}\label{eq. Dini}
h(x)\leq h(y) \implies h^\prime_-(y;y-x)\geq\frac{\gamma}{2}\|y-x\|^2,
\end{equation}
where $h^\prime_-(y;v):=\liminf\limits_{t\downarrow0}\frac{h(y+tv)-h(y)}{t}$ is the lower Dini derivative of $h$ at $y$ in the direction $v$, and that the converse holds with modulus $\frac{\gamma}{2}$. The same factor $\frac{1}{2}$ appears in \cite[Theorem 4 and Corollary 8]{hadjisavvas2026characterizations}, for a condition involving only values of $h$ and for G\^ateaux differentiable functions. Bao and Nam \cite[Section 5]{bao2025first} asked whether the modulus $\gamma$ can be recovered in this setting. The loss of the factor $\frac{1}{2}$ has a quantitative effect, since quadratic growth constants, error bounds and convergence rates depend on the modulus, and obtaining it only up to the factor $\frac{1}{2}$ leads to less sharp guarantees.

In another direction, the survey \cite{grad2025strongly} collects results on the Clarke subdifferential of these functions, obtained by Vial \cite[Lemma 1 and Theorems 2 and 3]{vial1982strong} and by Dudov and Osiptsev \cite[Lemma 2.1 and Corollaries 2.1 and 3.1]{dudov2022sufficient}, and records two open questions: one concerning the converse of Vial's inequality for Clarke subgradients \cite[Remark 49]{grad2025strongly}, and another concerning the converse of an inclusion between proximal points and strong subdifferentials \cite[Remark 38]{grad2025strongly}.

In this work, we answer these three questions. The main contribution is Theorem \ref{teo. P}, according to which, for $h$ continuous along segments of a convex set $C$, not necessarily open, condition \eqref{eq. Dini} with constant $\gamma$ is equivalent to strong quasiconvexity with modulus $\gamma$, even when the lower Dini derivative is replaced by the upper one. Thus, the exact modulus persists in the nonsmooth setting, and the converses mentioned above hold with the modulus $\gamma$ itself. The proof relies on a one-dimensional lemma (Lemma \ref{Lem. 1D}), which extends \cite[Lemma 3.1]{bao2025first} to functions that are merely continuous, and the result of Bao and Nam is recovered as a particular case (Corollary \ref{MR cor1}). As an application, we obtain a formula for the exact modulus of radial functions on Hilbert spaces (Corollary \ref{Coro radial}), which shows that the moduli of $\|x\|^\alpha$ that can be extracted from \cite{nam2026strong} are exact (Example \ref{Ex. Coro radial}), and we verify that the quadratic growth constant $\frac{\gamma}{4}$ of \cite[Corollary 9]{hadjisavvas2026characterizations} is optimal (Remark \ref{MR. obs2}).

The second contribution concerns the Clarke subdifferential. For locally Lipschitz functions on open convex subsets of Banach spaces, we prove that strong quasiconvexity with modulus $\gamma$ is equivalent to Vial's inequality with the same modulus, that is, to the inclusion of the Clarke subdifferential in the strong sublevel subdifferential $\partial_{1,\gamma}h$ of \cite{kabgani2022strong} (Theorem \ref{teo. Clarke}). This answers the question on Vial's inequality in the affirmative and extends \cite[Theorem 2]{vial1982strong} from $\mathbb{R}^n$ to Banach spaces. In $\mathbb{R}^n$, the same holds with the limiting subdifferential (Corollary \ref{Coro. subdlim}). It also follows that the classes of strongly quasiconvex functions of order $2$ considered in \cite{arora2011characterizations,singh2018strong} coincide with the class defined by \eqref{eq. SQ} (Remark \ref{rem. order}).

The third contribution concerns optimization. The direct implication of Theorem \ref{teo. Clarke} yields an optimality condition that generalizes \cite[Corollaries 2.1 and 3.1]{dudov2022sufficient} to functions that are not necessarily Clarke regular (Corollary \ref{Cons. coro_opt1}), an error bound with optimal constant (Corollary \ref{Cons. coro_opt2}), the rate $O(\frac{1}{k})$ for the squared distance to the minimizer in the projected subgradient method, which is the same as in the strongly convex case (Theorem \ref{Cons. teo1}), and the exponential convergence of the Clarke subgradient flow, which extends \cite[Theorem 15(b),(c)]{lara2025characterizations} to locally Lipschitz functions (Proposition \ref{Cons. Prop1}). The method uses only Clarke subgradients, which can be obtained through the usual calculus rules, and thus complements the method of \cite{choque2024subgradient}, which is based on strong subdifferentials. Moreover, Proposition \ref{Cons. Prop2} answers the question on proximal points, since it shows that the converse of \cite[Proposition 36]{kabgani2022strong} fails even for the strongly convex function $h(y)=y^2$, and that this converse holds for the truncation $\max\{h,h(x)\}$.

The paper is organized as follows. In Section \ref{Sec. Preli}, we gather the definitions and facts used on Dini derivatives and the Clarke subdifferential, including a monotonicity lemma. In Section \ref{Sec. MR}, we prove the one-dimensional lemma and, from it, the characterization of strong quasiconvexity via Dini derivatives, with its consequences for radial and differentiable functions, as well as the characterization via Clarke subgradients. In Section \ref{Sec. Conseq}, we derive an optimality condition, an error bound and the convergence of the projected subgradient method and of the associated flow, and we study the relation between proximal points and strong subdifferentials. Finally, Section \ref{Sec. Concl} presents the conclusions.

\section{Preliminaries}\label{Sec. Preli}

Throughout the text, $X$ denotes a real normed space, $X^*$ its topological dual, $\langle \cdot ,\cdot \rangle$ the pairing between $X^*$ and $X$, which coincides with the inner product when $X$ is a Hilbert space, and $B(x,r)$ the open ball with center $x$ and radius $r$. We now recall the definitions and auxiliary results used throughout the paper.

\subsection{Dini derivatives} \label{Subsec.Preli-Dini}
Let $C \subseteq X$ be convex and let $h: C \to \mathbb{R}$ be a function. For $x\in C$ and $a\in X$ such that $x+ta\in C$ for all sufficiently small $t>0$, the upper and lower Dini derivatives of $h$ at $x$ in the direction $a$ are given by
$$
h^\prime_+(x;a) = \limsup_{t\downarrow 0} \frac{h(x+ta)-h(x)}{t},
\quad \quad
h^\prime_-(x;a) = \liminf_{t\downarrow 0} \frac{h(x+ta)-h(x)}{t}.
$$
For a function $\psi: I \to \mathbb{R}$ defined on an interval $I \subseteq \mathbb{R}$, we use the same notation $\psi^\prime_{\pm}(s;r)$ for $s \in I$ and $r\in \mathbb{R}$ such that $s+tr\in I$ for all sufficiently small $t>0$, and we write $\psi^\prime_{\pm}(s):=\psi^\prime_{\pm}(s;1)$, that is,
\[
\psi^\prime_+(s) = \limsup_{t\downarrow 0} \frac{\psi(s+t)-\psi(s)}{t},
\quad \quad
\psi^\prime_-(s) = \liminf_{t\downarrow 0} \frac{\psi(s+t)- \psi(s)}{t}.
\]
These derivatives may take the values $\pm\infty$ and satisfy the following properties, which follow directly from the definitions.
\begin{itemize}
    \item [(P1)] Positive homogeneity, that is, $$h^\prime_\pm(x;\lambda a)=\lambda h^\prime_\pm(x;a)\quad \text{for all } \lambda>0.$$
    In particular, $\psi^\prime_\pm(s;r)=|r|\,\psi^\prime_\pm(s;\operatorname{sgn} r )$ for $r\neq 0$.
    \item [(P2)] $h^\prime_-(x; a)\leq h^\prime_+(x;  a)$.
    \item [(P3)] Let $x_0\in C$, $d\in X$ and $\psi(s)=h(x_0+sd)$ for $s$ in the interval $I:=\{s\in \mathbb{R}: x_0+sd\in C\}$, which is an interval containing $0$, since $C$ is convex. Then $\psi^\prime_\pm(s;r)=h^\prime_\pm(x_0+sd;rd)$ whenever one of the two sides is defined. In particular, with $X=\mathbb{R}$ and $d=-1$, the function $\chi(\sigma):=\psi(-\sigma)$ satisfies $\chi^\prime_\pm(-s)=\psi^\prime_\pm(s;-1)$.
\end{itemize}
We say that $h$ is continuous along segments of $C$ if, for all $x,y\in C$, the function $s\mapsto h(x+s(y-x))$ is continuous on $[0,1]$. This is the hypothesis of \cite[Theorems 4 and 6]{hadjisavvas2026characterizations}.

The following result is part of Dini's theorem on the Dini derivatives of continuous functions \cite[Chapter VI, p. 204]{saks1964theory}. Since it is short, we include the proof, which follows Saks' argument.

\begin{lemma}\label{Pre. Lem1}
Let $\psi:[a,b]\to \mathbb{R}$ be continuous and let $m\in\mathbb{R}$ be such that
$\psi^\prime_+(s)\geq m$ for all $s\in (a,b)$. Then
$$ \psi(b)-\psi(a)\geq m(b-a).$$
\end{lemma}
\begin{proof}
Let $m^\prime\in \mathbb{R}$ with $m^\prime<m$, $\varepsilon\in (0,b-a)$, and let $\varphi: [a,b]\to \mathbb{R}$ be given by $\varphi(s)=\psi(s)-m^\prime s$. Then $\varphi$ is continuous on $[a,b]$ and
$$ \varphi^\prime_+(s)=\psi^\prime_+(s)-m^\prime\geq m-m^\prime>0 \quad \text{for all } s\in(a,b).$$
Consider the set $S_\varepsilon:=\{s\in[a+\varepsilon,b]: \varphi(s)\geq \varphi(a+\varepsilon)\}$. Since $a+\varepsilon\in S_\varepsilon$, we have $S_\varepsilon\neq \emptyset$, and $S_\varepsilon$ is compact, since $\varphi$ is continuous. Let $M_\varepsilon:=\max S_\varepsilon$. If $M_\varepsilon<b$, then $M_\varepsilon\in (a,b)$, and the inequality $\varphi^\prime_+(M_\varepsilon)>0$ yields $s\in (M_\varepsilon,b]$ such that
$$\varphi (s)>\varphi(M_\varepsilon)\geq \varphi(a+\varepsilon),$$
that is, $s\in S_\varepsilon$ and $s>M_\varepsilon$, which contradicts the maximality of $M_\varepsilon$. Hence $M_\varepsilon=b$. Thus, $\psi(b)-\psi(a+\varepsilon)\geq  m^\prime(b-a-\varepsilon)$. Letting $\varepsilon \downarrow 0$, by the continuity of $\psi$ at $a$, and then $m^\prime \uparrow m$, we obtain
 $$\psi(b)-\psi(a)\geq  m (b-a).$$
\end{proof}

\begin{remark}
    Since $\psi^\prime_-\leq \psi^\prime_+$, Lemma \ref{Pre. Lem1} holds, in particular, under the hypothesis $\psi^\prime_-(s)\geq m$ for all $s\in(a,b)$.
\end{remark}

\subsection{Clarke subdifferential}\label{subs. clarke-prelim}

Let $X$ be a Banach space, $U\subseteq X$ open, $f:U\to \mathbb{R}$ and $x\in U$, with $f$ Lipschitz near $x$. The generalized directional derivative of $f$ at $x$ in the direction $v\in X$, denoted by $f^\circ (x;v)$, is defined by
$$f^\circ(x;v)=\limsup_{x^\prime \to x,~t\downarrow 0}\frac{f(x^\prime +tv)-f(x^\prime)}{t},$$
and the Clarke subdifferential of $f$ at $x$, denoted by $\partial^\mathcal{C}f(x)$, is the subset of $X^*$ given by
$$\partial^\mathcal{C}f(x):=\{\xi\in X^*:~\langle \xi,v\rangle\leq f^\circ(x;v) \quad \text{for all } v\in X \}.$$
These notions satisfy the following properties.
\begin{itemize}
\item [(P4)] The function $v\mapsto f^\circ(x;v)$ is finite and positively homogeneous, that is, $f^\circ(x;\lambda v)=\lambda f^\circ(x;v)$ for all $v\in X$ and $\lambda>0$.
\item [(P5)] $\partial^\mathcal{C}f(x)$ is a nonempty, convex and weak$^*$ compact subset of $X^*$.
\item [(P6)] $ f^\circ(x;v)=\max\{\langle \xi,v\rangle:~\xi \in \partial^\mathcal{C}f(x)\}$ for all $v \in X$.
\item [(P7)] $-f^\circ(x;-v)\leq f_-^\prime(x;v)\leq f_+^\prime(x;v)\leq f^\circ(x;v)$ for all $v \in X$.
\end{itemize}
Properties (P5)--(P6) can be found in \cite[Proposition 2.1.2]{clarke1990optimization}. In (P7), the middle inequality is (P2), and the last one follows from the definition of $f^\circ$ by taking $x^\prime=x$. For the first one, let $x_t:=x+tv$, so that $x_t \to x$ as $t \downarrow 0$. Since $x=x_t+t(-v)$,
 $$-f^\prime_-(x;v)=\limsup_{t \downarrow 0}\frac{f(x)-f(x+tv)}{t}=\limsup_{t \downarrow 0}\frac{f(x_t+t(-v))-f(x_t)}{t}\leq f^\circ(x;-v),$$
that is, $-f^\circ(x;-v)\leq f^\prime_-(x;v)$.

\section{Main results}\label{Sec. MR}
In this section, we prove a one-dimensional lemma and, from it, the characterization of strong quasiconvexity via Dini derivatives and via Clarke subgradients, without any loss in the modulus.

\subsection{A one-dimensional lemma}\label{subsec. Lem2} The following lemma will be used in the proof of Theorem \ref{teo. P}.

\begin{lemma}\label{Lem. 1D}
Let $\psi: [0,1]\to \mathbb{R}$ be continuous and $c>0$. Suppose that, for all $x\in (0,1)$ and $s\in [0,1]$ with $s\neq x$,
\begin{equation}\label{eq1 Lem1D}
\psi(s)\leq \psi(x)\implies \psi^\prime_+(x;x-s)\geq c(x-s)^2.
\end{equation}
Then
$$ \psi(t)\leq \max\{\psi(0),\psi(1)\}-ct(1-t)\quad \text{for all } t\in [0,1].$$
\end{lemma}
\begin{proof}
Let $M:=\max\{\psi(0),\psi(1)\}$, $c^\prime \in (0,c)$, and
$$\varphi(x):=\psi(x)+c^\prime x(1-x), \quad \quad x\in [0,1].$$
It suffices to show that $\varphi(t)\leq M$ for all $t\in [0,1]$ and then let $c^\prime \uparrow c$. Suppose, by contradiction, that $\max_{t\in [0,1]}\varphi(t)>M$. Since $\varphi$ is continuous, there exists $t^*\in[0,1]$ such that $\varphi(t^*)=\max_{t\in [0,1]}\varphi(t)>M$. Since $\varphi(0), \varphi(1)\leq M$, we have $t^*\in (0,1)$, and, for $a\in \{-1,1\}$ and all sufficiently small $t>0$, it follows that $t^*+ta\in[0,1]$ and $\varphi(t^*+ta)\leq \varphi(t^*)$. Hence,
$$  \limsup\limits_{t \downarrow 0}\frac{\varphi(t^*+ta)-\varphi (t^*)}{t}=\varphi^\prime_+(t^*;a)\leq 0.$$
Since $x\mapsto c^\prime x(1-x)$ is differentiable and $\varphi^\prime_+(t^*;a)=\psi^\prime_+(t^*;a)+c^\prime(1-2t^*) a$,
we obtain
    \begin{equation}\label{eq2 Lem1D}
        \psi^\prime_+(t^*;a)\leq-c^\prime (1-2t^*)\,a, \quad \quad a=\pm1.
    \end{equation}

\emph{Step 1.} We show that $t^*=\frac{1}{2}$. Suppose, by contradiction, that $t^*\neq\frac{1}{2}$, and let $s:=1-t^*\in(0,1)$ and $a:=t^*-s=2t^*-1\neq 0$. Since $s(1-s)=t^*(1-t^*)$ and $\varphi(s)\leq \varphi(t^*)$, it follows that $\psi(s)\leq\psi(t^*)$. Thus, \eqref{eq1 Lem1D} with $x=t^*$ guarantees that
$$\psi^\prime_+(t^*;a)\geq ca^2.$$
On the other hand, by property (P1) and \eqref{eq2 Lem1D}, we have
$$ ca^2\leq\psi^\prime_+(t^*;a)=|a|\,\psi^\prime_+(t^*;\operatorname{sgn} a)\leq-|a|\,c^\prime (1-2t^*)\operatorname{sgn} a=  c^\prime a^2.$$
Since $a\neq 0$, we obtain $c \leq c^\prime$, which contradicts the choice of $c^\prime$.

\emph{Step 2.} We describe the shape of $\psi$. Let $L_0:=\psi(\frac{1}{2})$. By \eqref{eq2 Lem1D} with $t^*=\frac{1}{2}$, we have $\psi^\prime_+(\frac{1}{2};\pm 1)\leq 0$. If $\psi(s)\leq L_0$ for some $s\neq \frac{1}{2}$, then \eqref{eq1 Lem1D} with $x=\frac{1}{2}$ and (P1) would imply
$$0<c\left(\tfrac{1}{2} -s\right)^2\leq\psi^\prime_+\left(\tfrac{1}{2};\tfrac{1}{2}-s\right)=  \left| \tfrac{1}{2}-s \right| \psi^\prime_+\left(\tfrac{1}{2};\operatorname{sgn}(\tfrac{1}{2}-s) \right)\leq 0,$$
which is absurd. Therefore,
\begin{equation}\label{eq3 Lem1D}
    \psi(s)>L_0\quad \text{for all } s\in[0,1]\setminus \{\tfrac{1}{2}\}.
\end{equation}
Moreover, since $\tfrac{1}{2}$ maximizes $\varphi$, we have $\varphi(x)\leq\varphi(\frac{1}{2})=L_0+\frac{c^\prime}{4}$ for all $x\in[0,1]$, that is,
\begin{equation}\label{eq4 Lem1D}
    \psi(x)\leq L_0+c^\prime\left(x-\tfrac{1}{2}\right)^2, \quad x\in [0,1].
\end{equation}
We now show that $\psi$ is strictly increasing on $[\frac{1}{2},1]$. For $x\in(\frac{1}{2},1)$, the inequality $\psi(\frac{1}{2})\leq \psi(x)$ allows us to apply \eqref{eq1 Lem1D} with $s=\frac{1}{2}$, and (P1) yields $\psi^\prime_+(x)\geq c (x-\tfrac{1}{2})$.
Let $\frac{1}{2}<\alpha<\beta\leq1$. Since $\psi^\prime_+\geq  c\left( \alpha-\tfrac{1}{2}\right)$ on $(\alpha, \beta)$, Lemma \ref{Pre. Lem1} guarantees that
 \begin{equation} \label{ast-01}\psi(\beta)-\psi(\alpha)\geq c\left(\alpha-\tfrac{1}{2}\right)(\beta-\alpha)>0,
\end{equation}
and, by \eqref{eq3 Lem1D}, $\psi(\beta)>\psi(\frac{1}{2})$ for all $\beta\in (\frac{1}{2},1]$, which proves the claim. The function $\chi(x):=\psi(1-x)$ is continuous, satisfies $\chi(\frac{1}{2})=L_0$, and fulfills \eqref{eq1 Lem1D} with the same constant $c$, since $\chi^\prime_+(x;a)=\psi^\prime_+(1-x;-a)$ by (P3). Applying the previous argument to $\chi$, we conclude that $\psi$ is strictly decreasing on $[0,\frac{1}{2}]$.

\emph{Step 3.} We estimate the width of the sublevel sets of $\psi$. Let $L_1:=\min\{\psi(0),\psi(1)\}$, which satisfies $L_1>L_0$ by \eqref{eq3 Lem1D}. By Step 2 and the intermediate value theorem, for each $L\in[L_0,L_1]$ there exist unique $x_-(L)\in[0,\frac{1}{2}]$ and $x_+(L)\in[\frac{1}{2},1]$ such that
\begin{equation}\label{eq4.1 Lem1D}
    \psi( x_-(L))=\psi( x_+(L))=L.
\end{equation}
Being inverses of the restrictions of $\psi$ to $[0,\frac{1}{2}]$ and to $[\frac{1}{2},1]$, the functions $x_-$ and $x_+$ are continuous, $x_-$ is decreasing, $x_+$ is increasing, and $x_\pm(L_0)=\frac{1}{2}$. Consequently, the function
 \begin{equation}\label{eq.w Lem1D}w(L):=x_+(L)-x_-(L)
\end{equation}
is continuous and increasing on $[L_0,L_1]$, with $w(L_0)=0$ and $w(L)>0$ for $L>L_0$. Applying \eqref{eq4 Lem1D} at $x=x_\pm(L)$, we obtain $|x_\pm(L)-\frac{1}{2}| \geq \sqrt{\frac{L-L_0}{c^\prime}}$, and adding these two bounds yields
\begin{equation}\label{eq5 Lem1D}
    w(L)\geq 2\sqrt{\frac{L-L_0}{c^\prime}}, \quad L\in [L_0,L_1].
\end{equation}

We now estimate the growth of $w$. Let $L_0<\eta <\eta^\prime \leq L_1$ and $\theta \in[x_+(\eta),x_+(\eta^\prime))$. Then $\theta\in (\frac{1}{2},1)$, $\psi(\theta)\in[\eta,\eta^\prime)$, and, by the uniqueness of $x_+$ in \eqref{eq4.1 Lem1D}, $\theta=x_+(\psi(\theta))$. Setting
$s:=x_-(\psi(\theta))$, we have $\psi(s)=\psi(\theta)$ and, by \eqref{eq.w Lem1D}, $\theta -s=w(\psi(\theta))>0$. Hence, by \eqref{eq1 Lem1D} with $x=\theta$, property (P1), and the monotonicity of $w$, we obtain
$$\psi^\prime_+(\theta)\geq c\,w(\psi(\theta))\geq c\,w(\eta), \quad  \theta \in [x_+(\eta),x_+(\eta^\prime)).$$
By Lemma \ref{Pre. Lem1} on $[x_+(\eta),x_+(\eta^\prime)]$, it follows that
\begin{equation}\label{eq6 Lem1D}
    \eta^\prime-\eta\geq  c\,w(\eta)\left(x_+(\eta^\prime)-x_+(\eta) \right).
\end{equation}
Similarly, if $\theta\in (x_-(\eta^\prime) ,x_-(\eta)]\subseteq (0,\frac{1}{2})$, then \eqref{eq1 Lem1D} with $s:=x_+(\psi(\theta))$ and (P1) imply $\psi^\prime_+(\theta;-1)\geq  c\,w(\eta)$. Since the function $\tilde\psi(\sigma):=\psi(-\sigma)$, $\sigma\in[-1,0]$, satisfies $\tilde \psi^\prime_+(-\theta) =\psi^\prime_+(\theta;-1)$ by (P3), Lemma \ref{Pre. Lem1} applied to $\tilde\psi$ on $[-x_-(\eta),-x_-(\eta^\prime)]$ yields
\begin{equation}\label{eq7 Lem1D}
    \eta^\prime-\eta\geq c\,w(\eta)\left ( x_-(\eta) -x_-(\eta^\prime)\right).
\end{equation}
Adding \eqref{eq6 Lem1D} and \eqref{eq7 Lem1D}, and recalling that $w$ is increasing, we have
\begin{equation}\label{eq8 Lem1D}
0\leq w(\eta^\prime)-w(\eta)\leq \frac{2(\eta^\prime-\eta)}{c\,w(\eta)}, \quad L_0<\eta<\eta^\prime \leq L_1.
\end{equation}

\emph{Step 4.} We conclude the proof. Fix $L_0<\eta_0<L\leq L_1$ and a partition $\eta_0<\eta_1<\cdots<\eta_N=L$ with $\Delta_k:=\eta_{k+1}-\eta_k\leq\delta$, and define $w_k:=w(\eta_k)$. Since $w_k\geq w_0>0$, it follows from \eqref{eq8 Lem1D} that
$$w_{k+1}^2-w_k^2=2w_k(w_{k+1}-w_k)+(w_{k+1}-w_k)^2\leq \frac{4\Delta_k}{c}+\frac{4\delta\, \Delta_k}{c^2w_0^2}.$$
Summing over $k=0,\ldots, N-1$ and using $\sum_k\Delta_k=L-\eta_0$, we obtain
$$w(L)^2-w(\eta_0)^2\leq \left(\frac{4}{c}+\frac{4\delta}{c^2w_0^2} \right) (L-\eta_0).$$
Since $w(L)$, $w(\eta_0)$ and $w_0$ do not depend on the partition, we may let $\delta\to 0$ in this inequality and then $\eta_0 \downarrow L_0$; the continuity of $w$ and $w(L_0)=0$ yield
 $w(L)^2\leq\frac{4}{c}(L-L_0)$. Hence, together with \eqref{eq5 Lem1D}, we obtain
$$\frac{4}{c^\prime}(L-L_0)\leq w(L)^2\leq \frac{4}{c}(L-L_0), \quad L\in(L_0,L_1],$$
and, since $L-L_0>0$, it follows that $c\leq c^\prime$, which contradicts the choice of $c^\prime$. Therefore, $\varphi(t)\leq M$ for all $t\in[0,1]$, which completes the proof.
 \end{proof}

Lemma \ref{Lem. 1D} extends \cite[Lemma 3.1]{bao2025first}, in which $\psi$ is differentiable on an open interval containing $[0,1]$ and the hypothesis is imposed with $\psi^\prime$ at all points of $[0,1]$; in our lemma, continuity suffices, and the hypothesis involves only the upper Dini derivative at interior points.

\begin{remark}\label{obs1}
\begin{enumerate}
\item [(i)] The constant $c$ in Lemma \ref{Lem. 1D} is optimal. Indeed, for $a>0$, the function $\psi(x)=a(x-\frac{1}{2})^2$ satisfies \eqref{eq1 Lem1D} if and only if $c\leq a$. With $u:=x-\frac{1}{2}$ and $v:=s-\frac{1}{2}$, the condition $\psi(s)\leq \psi(x)$ is equivalent to $|v|\leq |u|$, and
$$\psi^\prime_+(x;x-s)=2a\,u(u-v)=a(u-v)^2+a(u^2-v^2)\geq a(u-v)^2,$$
with equality when $s=1-x$. Moreover, this function satisfies $\psi(t)=\max\{\psi(0),\psi(1)\}-a\,t(1-t)$ for all $t\in [0,1]$.
\item [(ii)] The key point of the proof is Step 1, in which comparing $t^*$ with the symmetric point $1-t^*$ reduces the problem to the case where the maximizer of $\varphi$ is $\frac{1}{2}$. Bao and Nam \cite[Remark 3.3]{bao2025first} observe that controlling only the sum of the rates of change of the endpoints of the sublevel intervals, as in \eqref{eq8 Lem1D}, generally leads to half of the modulus; for this reason, in the differentiable case, they resort to a potential function that controls each endpoint separately. In Lemma \ref{Lem. 1D}, the reduction in Step 1 allows us to work with the bound on the sum. With the maximizer at $\frac{1}{2}$, the upper bound \eqref{eq4 Lem1D} is symmetric with respect to $\frac{1}{2}$, and the estimates for the width obtained in Steps 3 and 4,
$$\frac{4}{c^\prime}(L-L_0)\leq w(L)^2\leq \frac{4}{c}(L-L_0),$$
differ only by the replacement of $c^\prime$ by $c$.

\item [(iii)] When one of the endpoints of the segments is fixed at a minimizer, as in the case of strongly star quasiconvex functions \cite{khanh2025star}, studied under the name of strongly quasar-convex functions in \cite{nguyen2026projected} (see also \cite[Remark 13]{khanh2025star}), the differentiable characterization with the exact modulus is obtained in \cite[Theorem 3.3]{nguyen2026projected} and
 \cite[Proposition 27]{khanh2025star} by means of the mean value theorem or of an integration along the segment. In \eqref{eq. SQ}, both endpoints are arbitrary, and this is the situation handled by Lemma \ref{Lem. 1D}.
\end{enumerate}
\end{remark}

\subsection{Characterization via Dini derivatives}

Let $h:C\to \mathbb{R}$, with $C$ convex. The derivative $h^\prime_\pm(y;y-x)$ depends on the values of $h$ at the points $y+t(y-x)$ with small $t>0$, which may not belong to $C$ if $C$ is not open. We consider the set
$$\mathcal A(C):=\left \{(x,y)\in C\times C :~\text{there exists }\varepsilon>0 \text{ such that }y+t(y-x)\in C\text{ for all }t\in (0, \varepsilon)\right\}.$$
When $C$ is open, $\mathcal A(C)=C\times C$.

\begin{theorem}\label{teo. P}
Let $C$ be a convex subset of a real normed space, $h:C\to \mathbb R$ continuous along segments of $C$, and $\gamma >0$. Then the following statements are equivalent:
   \begin{enumerate}
       \item [(a)] $h$ is strongly quasiconvex with modulus $\gamma$;
       \item [(b)] for all $(x,y)\in  \mathcal{A}(C)$, implication \eqref{eq. Dini} holds;
       \item [(c)] for all $(x,y)\in \mathcal{A}(C)$, the implication
       \begin{equation}\label{eq. dini2}
           h(x)\leq h(y)  \Longrightarrow h^\prime_+(y;y-x)\geq \frac{\gamma}{2}\|y-x\|^2
       \end{equation}
       holds.
   \end{enumerate}
Moreover, the implication (a) $\Rightarrow$ (b) does not require the continuity of $h$.
\end{theorem}
\begin{proof}(a) $\Rightarrow$ (b). For $C$ open, this implication is \cite[Theorem 6(b)]{hadjisavvas2026characterizations}. We present a direct argument, which does not use the continuity of $h$. Let $(x,y)\in \mathcal{A}(C)$ with $h(x)\leq h(y)$, and let $\varepsilon>0$ be as in the definition of $\mathcal{A}(C)$. If $x=y$, both sides of \eqref{eq. Dini} vanish, so we may assume that $x\neq y$. Given $t\in (0,\varepsilon)$, the point $z_t:=y+t(y-x)$ belongs to $C$ and
    $$ y=\lambda z_t+(1-\lambda)x, \quad \lambda:=\frac{1}{1+t}, \quad  \|z_t-x\|=(1+t)\|y-x\|.$$
Since $\lambda(1-\lambda)=\frac{t}{(1+t)^2}$, inequality \eqref{eq. SQ} entails
$$h(y)\leq \max \{h(z_t),h(x)\}-\frac{\gamma}{2}\, t\|y-x\|^2.$$
If the maximum were $h(x)$, we would have $h(y)<h(x)$, which is incompatible with $h(x)\leq h(y)$. Therefore,
$$\frac{h(y+t(y-x))-h(y)}{t}\geq \frac{\gamma}{2}\|y-x\|^2 \quad \text{for all } t\in (0,\varepsilon),$$
and taking the lower limit as $t \downarrow 0$ yields \eqref{eq. Dini}.

\medskip

(b) $\Rightarrow$ (c). This follows from property (P2).

\medskip

(c) $\Rightarrow$ (a). Fix distinct $x_0,y_0\in C$, let $d:=y_0-x_0$ and $\psi(\tau):=h(x_0+\tau d)$ for $\tau\in[0,1]$. The function $\psi$ is continuous, since $h$ is continuous along segments. We verify that $\psi$ satisfies \eqref{eq1 Lem1D} with $c:=\frac{\gamma}{2}\|d\|^2>0$. Let $\tau\in (0,1)$ and $s\in [0,1]$ with $s\neq \tau$ and $\psi(s)\leq \psi(\tau)$, and consider $u:=x_0+\tau d$ and $v:=x_0+sd$. Then $h(v)\leq h(u)$ and, since $\tau \in (0,1)$, we have $u+t(u-v)=x_0+\left( \tau+t(\tau-s)\right)d \in[x_0,y_0]\subseteq C$ for all sufficiently small $t>0$, so that $(v,u)\in \mathcal A(C)$. Applying \eqref{eq. dini2} to the pair $(x,y)=(v,u)$ and using (P3), we obtain
 $$\psi^\prime_+(\tau;\tau-s)=h^\prime_+(u;u-v)\geq \frac{\gamma}{2}\|u-v\|^2=c\,(\tau-s)^2.$$
Thus, Lemma \ref{Lem. 1D} guarantees that
$$h\left( (1-t)x_0+ty_0\right)\leq \max\{h(x_0), h(y_0)\}-\frac{\gamma}{2}\, t(1-t)\|y_0-x_0\|^2, \quad t\in [0,1],$$
which is \eqref{eq. SQ}. Since $x_0$ and $y_0$ are arbitrary, and the case $x_0=y_0$ is trivial, $h$ is strongly quasiconvex with modulus $\gamma$.
\end{proof}

Theorem \ref{teo. P} holds in an arbitrary real normed space, as does \cite[Theorem 6]{hadjisavvas2026characterizations}, without requiring the domain $C$ to be open. Moreover, it shows that the best constant in \eqref{eq. Dini} coincides with the largest modulus in \eqref{eq. SQ}, whereas the converse of \cite[Theorem 6(b)]{hadjisavvas2026characterizations} only ensures the modulus $\frac{\gamma}{2}$.

It is already observed in \cite{vladimirov1978uniformly}, right after Theorem 6, that it suffices to assume absolute continuity along segments and to impose the inequality at the points where the directional derivative exists, with modulus $\frac{\gamma}{2}$. Theorem \ref{teo. P}, in turn, requires only continuity along segments, uses Dini derivatives, and ensures the modulus $\gamma$.

Theorem \ref{teo. P} also allows one to compute the exact modulus from first-order information. The following result deals with radial functions.

\begin{corollary}\label{Coro radial}
Let $H\neq \{0\}$ be a real Hilbert space, $R>0$, $C:=\{x \in H: \|x\|\leq R\}$, $f:[0,R]\to \mathbb{R}$ continuous and strictly increasing, and let $h:C\to \mathbb R$ be the radial function $h(x):=f(\|x\|)$. Let
\begin{equation}\label{eq1 Coro par}
    \gamma^*:=\inf_{0<r<R}\frac{f^\prime_-(r)}{r},
\end{equation}
where $f^\prime_-$ is the lower Dini derivative of $f$.
Then, for $\gamma>0$, the function $h$ is strongly quasiconvex with modulus $\gamma$ if and only if $\gamma \leq \gamma^*$.
\end{corollary}
\begin{proof}
Since $C$ is convex and $h$ is continuous, by Theorem \ref{teo. P}
it suffices to show that \eqref{eq. Dini} holds for all $(x,y)\in \mathcal{A}(C)$ if and only if $\gamma\leq \gamma^*$.
Let $(x,y)\in \mathcal{A}(C)$ with $x\neq y$ and $\|x\|\leq \|y\|$, a condition equivalent to $h(x)\leq h(y)$, and let $d:=y-x$. Then $0<\|y\|<R$, since $y=0$ would imply $x=y$, and, if $\|y\|=R$, the Cauchy--Schwarz inequality would give $\langle y,d\rangle>0$, hence $\|y+td\|>R$ for all $t>0$.
Let $r:=\|y\|$, $u:=\frac{y}{r}$, $a:=\langle x, u\rangle$, $b:=x-au$, which satisfies $b\perp u$, and $c:=\langle u, d\rangle = r-a$. Then $c>0$, since $c=0$ would imply $x=y$, and
$$\|b\|^2=\|x\|^2-a^2 \leq c(2r-c).$$
From $\|y+td\|^2=(r+tc)^2+t^2\|b\|^2$ it follows that, for each
$\varepsilon>0$ and all sufficiently small $t>0$,
 $ r+tc\leq \|y+td\|\leq r+t(c+\varepsilon)$.
Since $f$ is increasing, it follows that $ h^\prime_-(y;d)= c f^\prime_-(r)$. Since $\|d\|^2=c^2+\|b\|^2$, condition \eqref{eq. Dini} is equivalent, for this pair, to
\begin{equation}\label{eq1 Coro radial}
    f^\prime_-(r)\geq \frac{\gamma}{2}\left(c+\frac{\|b\|^2}{c} \right).
\end{equation}
Moreover, the inequality $\|b\|^2\leq c(2r-c)$ gives
$c+ \frac{\|b\|^2}{c}\leq 2r$, with equality when $\|x\|=\|y\|$, in particular when $x=-y$. It thus remains to prove that \eqref{eq1 Coro radial} holds for all $(x,y)\in \mathcal{A}(C)$ with $x\neq y$ and $\|x\|\leq\|y\|$ if and only if $\gamma\leq \gamma^*$.

($\Leftarrow$) If $\gamma \leq \gamma^*$, then, for all such pairs $(x,y)$, we have $$f^\prime_-(r)\geq \gamma r\geq \frac{\gamma}{2}\left(c+\frac{\|b\|^2}{c} \right),$$ that is, \eqref{eq1 Coro radial} holds.

($\Rightarrow$) Suppose that \eqref{eq. Dini} holds for all $(x,y)\in\mathcal{A}(C)$ with $x\neq y$, and let $r\in (0,R)$. Since $H\neq \{0\}$, we may take $y\in C$ with $\|y\|=r$. Then $(-y,y)\in  \mathcal{A}(C)$, with $c=2r$ and $b=0$, and \eqref{eq. Dini} reduces to $f^\prime_-(r)\geq  \gamma r$. Therefore, $\gamma\leq \gamma^*$.
\end{proof}

Corollary \ref{Coro radial} reduces the computation of the exact modulus of a radial function to the computation of an infimum in one variable. The following two examples illustrate this reduction. The first one deals with powers of the norm, which include nonconvex and non-Lipschitz functions, and the second one compares, for functions of one variable, the exact moduli with the moduli guaranteed by the converse with factor $\frac{1}{2}$.

\begin{example}\label{Ex. Coro radial}
Let $H\neq \{0\}$, $R>0$ and $C$ be as in Corollary \ref{Coro radial}. For $\alpha>0$, let $h(x)=f(\|x\|)$ with $f(r):=r^\alpha$. Since $f^\prime_-(r)=\alpha r^{\alpha -1}$ for $r>0$, it follows from Corollary \ref{Coro radial} that
$$ \gamma^*=\inf_{0<r<R} \alpha\, r^{\alpha-2}=
\begin{cases}
\alpha R^{\alpha-2}, & 0<\alpha<2,\\
2, & \alpha=2,\\
0, & \alpha>2 .
\end{cases}$$
Thus, for $0<\alpha\leq 2$, $h$ is strongly quasiconvex on $C$ with exact modulus $\gamma^*$, whereas, for $\alpha>2$, $h$ is not strongly quasiconvex on $C$, whatever $R>0$ is. For $0<\alpha<1$, $h$ is neither convex nor locally Lipschitz at $0$, and its strong quasiconvexity on bounded convex sets is known \cite[Example 8]{grad2025strongly}. The proofs of \cite[Proposition 3.1 and Corollary 3.9]{nam2026strong} lead to the modulus $\alpha R^{\alpha-2}$ for $0<\alpha\leq 1$, and Corollary \ref{Coro radial} shows that this value is exact. The case $1<\alpha\leq 2$ is mentioned in Remark 3.10 of the same work.
\end{example}

\begin{example}\label{exem. Table}
    Table \ref{Tab. const} compares the exact moduli with those guaranteed by \cite[Theorem 6(b)]{hadjisavvas2026characterizations}, applied to the interior of each interval, which does not change the moduli, since the functions are continuous. Functions (ii)--(iv) are even, and for them we use Corollary \ref{Coro radial} with $H=\mathbb{R}$. Functions (i) and (v) are not even, and for them we apply Theorem \ref{teo. P} directly. In these two cases, for $y>0$, the points $x\neq y$ with $h(x)\leq h(y)$ satisfy $x<y$, and the largest admissible constant in \eqref{eq. Dini} for the pair $(x,y)$, which equals $\frac{2h^\prime(y)}{y-x}$, is smallest when $x$ is the smallest of these points. Pairs with $y\leq 0$ produce larger bounds.
\begin{enumerate}
    \item [(i)] For $h(t)=t^2$ $(t\geq 0)$ and $h(t)=4t^2$ $(t<0)$ on $[-3,3]$ and $0<y<3$, the smallest of these points is $x=-\frac{y}{2}$. Thus,
    $$\gamma^* =\frac{2\cdot 2y}{y+\frac{y}{2}}=\frac{8}{3}.$$
    \item [(ii)] For $h(t)=|t|+t^2$ on $[-1,1]$, we have $f(r)=r+r^2$ and $\frac{f^\prime_-(r)}{r}=2+\frac{1}{r}$, whose infimum on $(0,1)$ is $\gamma^*=3$.
    \item [(iii)] For $h(t)=\max\{t^2,3t^2-1\}$ on $[-1,1]$, we have $\frac{f^\prime_-(r)}{r}=2$ for $0<r<\frac1{\sqrt 2}$ and $\frac{f^\prime_-(r)}{r}=6$ for $\frac{1}{\sqrt{2}}\leq r<1$. Hence, $\gamma^*=2$.
    \item [(iv)] For $h(t)=t^2+2\sin^2 t$ on $[-2,2]$, the function $f(r)=r^2+2\sin^2 r$ is strictly increasing on $[0,\infty)$, since $f^\prime(r)=2r+2\sin 2 r>0$ for $r>0$, and $\frac{f^\prime_-(r)}{r}=2+\frac{2\sin(2r)}{r}$. The function $u\mapsto \frac{\sin u}{u}$ is decreasing on $(0,u^*]$ and attains at $u^*$ its minimum on $(0,\infty)$, where $u^*\approx 4.4934$ is the first positive solution of $\tan u=u$. Since $2r<4<u^*$ for $0<r<2$, we obtain $\gamma^*=2+\sin 4\approx 1.2432$. Moreover, $h$ is strongly quasiconvex on $\mathbb{R}$ with modulus $\gamma$ if and only if it is so on each interval $[-R,R]$, and therefore its exact modulus on the whole real line is
    $$\gamma^*_{\mathbb{R}}=\inf_{r>0}\left(2+\frac{2\sin 2r}{r}\right)=2+4\cos u^*\approx 1.1311.$$
    In particular, the modulus $\frac{1}{2}$ considered in \cite[Remark~7(ii) and Example~30]{lara2025characterizations} is admissible, and the exact modulus on $\mathbb{R}$ is more than twice this value.
    \item [(v)] For $h(t)=\operatorname{sgn}(t)\sqrt{|t|}$ on $[-\frac{1}{2},\frac{3}{2}]$ and $0<y<\frac{3}{2}$, the smallest of these points is $x=-\frac{1}{2}$, and we obtain
    $$\gamma^*=\inf_{0<y<\frac{3}{2}}\frac{2h^\prime(y)}{y+\frac{1}{2}}=\lim_{y\uparrow \frac{3}{2}}\frac{1}{\sqrt y\,(y+\frac12)}=\frac{1}{\sqrt6}.$$
\end{enumerate}
\end{example}

\begin{table}[H] \centering \small \caption{Exact moduli of strong quasiconvexity and moduli guaranteed by \cite[Theorem~6(b)]{hadjisavvas2026characterizations}.}\label{Tab. const} \begin{tabular}{llcc} \toprule $h(t)$ & interval & exact & by \cite{hadjisavvas2026characterizations}\\ \midrule $t^2$ ($t\ge0$), $4t^2$ ($t<0$) & $[-3,3]$ & $8/3$ & $4/3$\\ $|t|+t^2$ & $[-1,1]$ & $3$ & $3/2$\\ $\max\{t^2,3t^2-1\}$ & $[-1,1]$ & $2$ & $1$\\ $t^2+2\sin^2t$ & $[-2,2]$ & $2+\sin4\approx1.2432$ & $0.6216$\\ $\operatorname{sgn}(t)\sqrt{|t|}$ & $[-\frac{1}{2},\frac{3}{2}]$ & $1/\sqrt{6}\approx0.4082$ & $0.2041$\\ \bottomrule \end{tabular} \end{table}

\begin{remark}\label{MR. obs2}
If $h$ satisfies \eqref{eq. SQ} and $\bar x$ minimizes $h$ on $C$, then
$$ h(y)-h(\bar x)\geq \frac{\gamma}{4}\|y-\bar x\|^2\quad \text{for all } y\in C $$
\cite[Corollary 9]{hadjisavvas2026characterizations}
(for $h$ lower semicontinuous, see also \cite[Theorem 5.1]{nam2026strong}). The factor $\frac{1}{4}$ is optimal. Indeed, let $R>0$, $b>1$, $C=[-R,R]$, and $h_b(t):=t^2$ if $t\geq 0$ and $h_b(t):=bt^2$ if $t<0$. Theorem \ref{teo. P}, applied to the extreme cases $x=-\frac{y}{\sqrt{b}}$, with $0<y<R$, and $x=-\sqrt{b}\, y$, with $-\frac{R}{\sqrt{b}}\leq y< 0$, shows that the exact modulus of $h_b$ is
$$\gamma_b=\min\left\{ \frac{4\sqrt{b}}{1+\sqrt{b}},~\frac{4b}{1+\sqrt{b}} \right\}=\frac{4\sqrt{b}}{1+\sqrt{b}},$$
since all other pairs produce larger bounds. For $b=4$ and $R=3$, one obtains $\gamma_4=\frac{8}{3}$, as in row (i) of Table~\ref{Tab. const}. On the other hand, $\bar x=0$ and
$$\min_{t\in C\setminus \{0\}}\frac{h_b(t)-h_b(0)}{t^2}=1=\frac{1+\frac{1}{\sqrt {b}}}{4}\,\gamma_b,
\quad  \quad
\frac{1+\frac{1}{\sqrt {b}}}{4}\ \longrightarrow\ \frac{1}{4}\quad(b\to\infty).$$
The asymmetry of $h_b$ is essential. For radial functions, under the assumptions of Corollary~\ref{Coro radial}, the inequality $f^\prime_-(r)\geq\gamma^*r$ for $0<r<R$ and Lemma \ref{Pre. Lem1}, applied to $r\mapsto f(r)-\frac{\gamma^*}{2}r^2$, imply $f(r)-f(0)\geq\frac{\gamma^*}{2}r^2$ for all $r\in[0,R]$, that is, $h(y)-h(0)\geq\frac{\gamma^*}{2}\|y\|^2$ for all $y\in C$, with the constant $\frac{\gamma^*}{2}$ in place of $\frac{\gamma^*}{4}$.
\end{remark}

Recall that $h$, defined on an open set $C\subseteq X$, is G\^ateaux differentiable at $y\in C$ if there exists $\nabla h(y)\in X^*$ such that $\lim_{t\to 0}\frac{h(y+tv)-h(y)}{t}=\langle{\nabla h(y)}, v\rangle$ for all $v\in X$.
For such functions, Theorem~\ref{teo. P} reduces to the following result.

\begin{corollary}\label{MR cor1}
Let $C$ be an open convex subset of a real normed space, $h:C\to \mathbb{R}$ G\^ateaux differentiable, and $\gamma>0$. Then $h$ is strongly quasiconvex with modulus $\gamma$ if and only if \eqref{eq. grad} holds for all $x,y\in C$.
\end{corollary}
\begin{proof} Given $x,y\in C$, since $C$ is open, the function $\tau \mapsto h(x+\tau(y-x))$ is differentiable on an open interval containing $[0,1]$, so that $h$ is continuous along segments of $C$, and
$$h^\prime_-(y;y-x)=\langle \nabla h(y),y-x\rangle.$$
Since $\mathcal{A}(C)=C\times C$, the result follows from the equivalence (a) $\Leftrightarrow$ (b) of Theorem \ref{teo. P}.
\end{proof}

Corollary \ref{MR cor1} coincides with \cite[Theorem 3.2]{bao2025first} in the case $\gamma>0$. In this work, it is obtained as a particular case of Theorem \ref{teo. P}, through a different argument,
based solely on Dini derivatives. Likewise, if $h$ is of class $C^1$ on an open set containing a convex set $C$, not necessarily open, then \eqref{eq. grad} implies \eqref{eq. Dini} on $\mathcal A(C)$, and Theorem \ref{teo. P} guarantees the modulus $\gamma$ in the converse of \cite[Theorem~6]{vladimirov1978uniformly}, as in the formulation of \cite[Proposition 40]{grad2025strongly}.

Analogously, the characterization in \cite[Theorem 4]{hadjisavvas2026characterizations}, which involves only values of $h$ at three collinear points $x$, $y$ and $z \in [x,y]$, and whose converse is established in that work with modulus $\frac{\gamma}{2}$, also holds with the same modulus in both directions.

\begin{corollary}\label{MR cor2}
    Let $C$ be a convex subset of a real normed space, $h:C\to \mathbb{R}$ continuous along segments of $C$, and $\gamma>0$. Then $h$ is strongly quasiconvex with modulus $\gamma$ if and only if, for all $x,y\in C$ and $z=x+t(y-x)$ with $0<t\leq 1$, the implication
    \begin{equation}\label{eq.1 Cor2}
h (x)\leq h(z)\quad \Longrightarrow \quad h(z)\leq h(y)-\frac{\gamma}{4}(1-t^2)\|y-x\|^2
\end{equation}
holds.
\end{corollary}

\begin{proof} ($\Rightarrow$) This implication is \cite[Theorem 4]{hadjisavvas2026characterizations}, which holds on any convex set and does not use the continuity of $h$.

($\Leftarrow$) By Theorem \ref{teo. P}, it suffices to verify item (b). Let $(x,y)\in \mathcal{A}(C)$ with $h(x)\leq h(y)$. If $x=y$, both sides of \eqref{eq. Dini} vanish. Suppose, then, that $x\neq y$, and let $y_s:=y+s(y-x)$, which belongs to $C$ for all sufficiently small $s>0$. Since $y=x+\frac{1}{1+s}(y_s-x)$ and $\|y_s-x\|=(1+s)\|y-x\|$, implication \eqref{eq.1 Cor2}, applied to the points $x$ and $y_s$ with $z=y$ and $t=\frac{1}{1+s}$, yields
$$ h(y)\leq h(y_s)-\frac{\gamma}{4}\left( \|y_s-x\|^2-\|y-x\|^2\right)=h(y_s)-\frac{\gamma}{4}\, (2s+s^2)\|y-x\|^2.$$
Dividing by $s$, we obtain
$$\frac{h(y_s)-h(y)}{s}\geq \frac{\gamma}{4}\,(2+s)\|y-x\|^2,$$
and taking the lower limit as $s\downarrow 0$ gives $h^\prime_-(y;y-x)\geq \frac{\gamma}{2}\|y-x\|^2$.
Thus, item (b) of Theorem \ref{teo. P} holds and, since $h$ is continuous along segments, $h$ is strongly quasiconvex with modulus $\gamma$.
\end{proof}

\subsection{Characterization via Clarke subgradients}\label{Sec. Clarke}

In this subsection, we answer the question in \cite[Remark 49]{grad2025strongly} by showing that Vial's inequality for Clarke subgradients characterizes strong quasiconvexity with the same modulus. Necessity follows directly from the definition of $h^\circ$, and sufficiency reduces to Theorem \ref{teo. P}.

Let $C$ be a subset of a real Banach space $X$ and $h:C\to \mathbb{R}$. Following \cite[Definition 19]{kabgani2022strong}, in the equivalent form of \cite[relation (7)]{lara2026optimality}, given $x\in C$,
$\beta>0$ and $\gamma\geq 0$, the strong sublevel subdifferential of $h$ at $x$ is
$$\partial_{\beta,\gamma}h(x):=\left\{ z\in X^*: \langle z,y-x\rangle \leq -\frac{\beta\gamma }{2}\|y-x\|^2\quad \text{for all } y\in S_x\right\},$$
where $S_x:=\{y\in C: h(y)\leq h(x)\}$. In \cite{kabgani2022strong,lara2026optimality}, proper functions $h:\mathbb{R}^n \to \mathbb{R}\cup\{+\infty\}$ and the sublevel set $\{y\in \mathbb{R}^n: h(y)\leq h(x)\}$ are considered, which coincides with $S_x$ when $h$ is extended by $+\infty$ outside $C$. In this case, as observed right after \cite[Definition 19]{kabgani2022strong}, $\partial_{\beta,\gamma}h(x)$ coincides with the strong subdifferential $\partial^{K}_{\beta,\gamma}h(x)$ of Definition 5 of the same work, with $K=S_x$, and the definition makes sense, without changes, in Banach spaces.

Being an intersection of weak$^*$ closed half-spaces, $\partial_{\beta,\gamma}h(x)$ is convex and weak$^*$ closed, and it depends on $h$ only through $S_x$. For $\gamma=0$, it is the normal cone to $S_x$ at $x$ in the sense of convex analysis. For $C$ open and $h$ locally Lipschitz, Vial's inequality can be written, in this notation, as the inclusion
 $$\partial^\mathcal{C}h(x)\subseteq \partial_{1,\gamma}h(x) \quad \text{for all } x\in C.$$

\begin{theorem}\label{teo. Clarke}
Let $C$ be an open convex subset of a real Banach space $X$, $h:C\to \mathbb{R}$ locally Lipschitz, and $\gamma>0$. The following statements are equivalent:
\begin{itemize}
\item [(a)] $h$ is strongly quasiconvex with modulus $\gamma$;

\item [(b)] for all $x,y\in C$, the implication $$h(y)\leq h(x)  \Longrightarrow h^\circ(x;y-x)\leq -\frac{\gamma}{2}\|y-x\|^2 $$ holds;

\item[(c)] $\partial^\mathcal{C}h(x)\subseteq \partial_{1,\gamma}h(x)$ for all $x\in C$, that is, for all $x,y\in C$,
$$ h(y)\leq h(x)\quad \Longrightarrow \quad \langle z,y-x\rangle\leq -\frac{\gamma}{2}\|y-x\|^2   \text{ for all } z\in \partial^\mathcal{C}h(x).$$
\end{itemize}
\end{theorem}
\begin{proof}
(a) $\Rightarrow$ (b). Let $x,y\in C$ with $h(y)\leq h(x)$ and $d:=y-x$. If $y=x$, the implication is immediate, since $h^\circ(x;0)=0$. Assume $y\ne x$ and consider two cases.

\emph{Case $h(y)<h(x)$.} Since $C$ is open and $h$ is locally Lipschitz, hence continuous, there exist $r>0$ and $\ell>0$ such that $B(x,r)\subseteq C$, $h$ is $\ell$-Lipschitz on $B(x,r)$, and $h(z)>h(y)$ for all $z\in B(x,r)$. If $x^\prime\in B(x,\frac{r}{2})$ and $t>0$ is sufficiently small, the points $x^\prime+td$ and $x^\prime+t(y-x^\prime)$ belong to $B(x,r)$, and they differ by $t(x^\prime-x)$. Since $h(y)<h(x^\prime)$, inequality \eqref{eq. SQ}, applied to the points $x^\prime$ and $y$, yields
$$\frac{h(x^\prime+td)-h(x^\prime)}{t} \leq \frac{h(x^\prime+t(y-x^\prime))-h(x^\prime)}{t}+\ell\|x^\prime-x\| \leq -\frac{\gamma}{2}(1-t)\|y-x^\prime\|^2+\ell\|x^\prime-x\| .$$
Taking the upper limit as $x^\prime\to x$ and $t\downarrow 0$, we obtain $h^\circ(x;d)\leq-\frac{\gamma}{2}\|d\|^2$.

\emph{Case $h(y)=h(x)$.} For $s\in(0,1)$, let $y_s:=x+sd$. By \eqref{eq. SQ},
$$h(y_s) \leq h(x)-\frac{\gamma}{2}\,s(1-s)\|d\|^2\ <\ h(x).$$
Applying the previous case to the pair $x$, $y_s$ and using (P4), we get
$$s\,h^\circ(x;d)=h^\circ(x;sd) \leq -\frac{\gamma}{2}\,s^2\|d\|^2 .$$
Dividing by $s$ and letting $s \uparrow 1$, we conclude that $h^\circ(x;d) \leq -\frac{\gamma}{2}\|d\|^2$.

(b) $\Rightarrow$ (c). By the definition of $\partial^{\mathcal C}h(x)$, we have $\langle z,y-x \rangle \leq  h^\circ(x;y-x)$ for all $z\in \partial^{\mathcal C}h(x)$, and (c) follows from (b).

(c) $\Rightarrow$ (a). Let $x,y\in C$ with $h(x)\leq h(y)$. Since $x\in S_y$, the inclusion $\partial^{\mathcal C}h(y)\subseteq\partial_{1,\gamma}h(y)$ yields
$$\langle z,{x-y}\rangle \leq-\frac{\gamma}{2}\|y-x\|^2 \quad  \text{for all } z\in\partial^{\mathcal C}h(y).$$
Hence, by (P6), $h^\circ(y;x-y)\leq-\frac{\gamma}{2}\|y-x\|^2$. By (P7) with $v=y-x$, we obtain
$$h^\prime_-(y;y-x)\ \ge\ -h^\circ(y;x-y) \geq \frac{\gamma}{2}\|y-x\|^2.$$
Since $C$ is open, we have $\mathcal A(C)=C\times C$, and $h$ is continuous, being locally Lipschitz. Therefore, item (b) of Theorem~\ref{teo. P} holds, which is equivalent to (a). This completes the proof.
\end{proof}

\begin{remark}
\begin{enumerate}
\item [(i)] The implication (a) $\Rightarrow$ (c) was obtained by Vial \cite[Theorem 2]{vial1982strong} for $C=X=\mathbb{R}^n$ and, for $C\subseteq\mathbb{R}^n$ open and convex and $h$ Clarke regular, it follows from \cite[Lemma 2.1]{dudov2022sufficient}. Analogous results for Fr\'echet subgradients and for the strong subdifferential can be found in \cite[Proposition 1]{lara2026optimality} and \cite[Proposition 40]{kabgani2022strong}, respectively.

\item [(ii)] In $\mathbb{R}^n$, the implication (a) $\Rightarrow$ (c) can also be obtained from \cite[Proposition 3(a)]{lara2026optimality}, applied with $V=C$, and from the equality $\partial^{\mathcal{C}}h=\conv \partial h$ \cite[Theorem 9.61]{rockafellar1998variational}, where $\conv$ denotes the convex hull. To this end, it suffices to verify that $y\mapsto S_y$ is inner semicontinuous at each $x\in C$, that is, $S_x\subseteq\liminf_{y\to x}S_y$. Let $x\in C$, $z\in S_x \setminus \{x\}$, and $z_s:=x+s(z-x)$ with $s\in(0,1)$. By \eqref{eq. SQ},
\[
h(z_s)\leq h(x)-\frac{\gamma}{2}\,s(1-s)\|z-x\|^2<h(x).
\]
Since $h$ is continuous, $z_s\in S_y$ for all $y$ sufficiently close to $x$, that is, $z_s\in\liminf_{y\to x}S_y$. Since this inner limit is closed and $z_s\to z$ as $s\uparrow1$, we have $z\in\liminf_{y\to x}S_y$, and the case $z=x$ follows from $y\in S_y$.

\item [(iii)] The implication (c) $\Rightarrow$ (a) is the converse of \cite[Proposition 48]{grad2025strongly}, which, by Theorem \ref{teo. Clarke}, holds with the same modulus and in Banach spaces.

\item [(iv)] For $h$ of class $C^1$, item (c) reduces to $\nabla h(x)\in\partial_{1,\gamma}h(x)$, and, for $C$ open, Theorem \ref{teo. Clarke} contains the characterization of \cite[Corollary~41]{grad2025strongly}.

\item [(v)] If $C$ is convex with nonempty interior, $h$ is continuous on $C$ and locally Lipschitz on $\intr C$, then $h$ is strongly quasiconvex with modulus $\gamma$ on $C$ if and only if $\partial^{\mathcal C}h(x)\subseteq\partial_{1,\gamma}h(x)$ for all $x\in\intr C$. Indeed, Theorem \ref{teo. Clarke} on $\intr C$ yields \eqref{eq. SQ} on $\intr C$, which extends to $C\subseteq\overline{\intr C}$ by the continuity of $h$. At boundary points, the inclusion may fail, since $\partial^{\mathcal C}h$ depends on the values of $h$ outside $C$. For example, let $C=[0,1]$, $h(t)=t^2$ if $t\le1$ and $h(t)=2-t$ if $t>1$. Then $h$ is strongly quasiconvex on $C$ and $\partial^{\mathcal C}h(1)=[-1,2]$. Since $0\in S_1$ and $\langle-1, 0-1\rangle=1>-\frac{\gamma}{2}$, we have $-1\notin\partial_{1,\gamma}h(1)$ for all $\gamma>0$.
\end{enumerate}
\end{remark}

In $\mathbb{R}^n$, Theorem \ref{teo. Clarke} admits a version with the limiting subdifferential, which we now recall \cite[Definition 8.3]{rockafellar1998variational}. Let $U\subseteq \mathbb{R}^n$ be open, $h:U\to \mathbb{R}$ locally Lipschitz, and $x\in U$. The regular subdifferential of $h$ at $x$ is
$$\hat\partial h(x):=\left\{\xi\in \mathbb{R}^n: h(y)\geq h(x)+\langle \xi, y-x \rangle +o(\|y-x\|)\right\},$$
and the limiting subdifferential of $h$ at $x$ is
$$\partial h(x):=\left\{\xi \in \mathbb{R}^n: \text{there exist } x_k \to x \text{ and } \xi_k \in \hat\partial h(x_k) \text{ with }\xi_k \to \xi \right\}.$$
Since $h$ is continuous, the condition $h(x_k)\to h(x)$ required in this definition is automatically satisfied. Moreover, by Theorem 9.61 of the same book, $\partial^{\mathcal C}h(x)=\conv \partial h(x)$.
\begin{corollary}\label{Coro. subdlim}
    Let $C\subseteq \mathbb{R}^n$ be open and convex, $h:C\to \mathbb{R}$ locally Lipschitz, and $\gamma>0$. Then $h$ is strongly quasiconvex with modulus $\gamma$ if and only if $ \partial h(x)\subseteq\partial_{1,\gamma}h(x)$ for all $x\in C$.
\end{corollary}
\begin{proof}
Since $\partial^{\mathcal{C}}h(x)=\conv \partial h(x)$ and $\partial_{1,\gamma}h(x)$ is convex, we have
$$\partial h(x)\subseteq\partial_{1,\gamma}h(x)\iff \partial^{\mathcal{C}}h(x)\subseteq\partial_{1,\gamma}h(x),$$
and it suffices to apply Theorem \ref{teo. Clarke}.
\end{proof}

\begin{remark}\label{rem. order}
Let $C\subseteq \mathbb{R}^n$ be open and convex. For a differentiable function $f:C\to \mathbb{R}$, \cite[Definition 14.4]{arora2011characterizations} calls $f$ strongly quasiconvex of type I of order $\sigma$, and
\cite[Definition 7]{singh2018strong} calls it strongly quasiconvex of order $\sigma$, if there exists $\alpha>0$ such that, for all $x,y\in C$ with $x\neq y$,
$$f(x)\leq f(y) \Longrightarrow \langle \nabla f(y), x-y\rangle +\alpha\|x -y\|^\sigma \leq 0.$$
For $f$ locally Lipschitz, \cite[Definition 5.1]{singh2021strong} requires, for all $x,y\in C$ and all $\eta\in\partial^{\mathcal{C}}f(y)$, some $\alpha=\alpha(x,y,\eta)>0$ such that
$$f(x)\leq f(y)  \Longrightarrow  \langle \eta, x-y\rangle +\alpha\|x-y\|^\sigma\leq 0.$$
For $\sigma=2$ and constant $\alpha$, Corollary \ref{MR cor1} and
Theorem \ref{teo. Clarke} show that these classes coincide with that of strongly quasiconvex functions with $\gamma=2\alpha$. With $\alpha$ depending on $(x,y,\eta)$, the class is larger. For example, with $C=\mathbb{R}$, the function $f(t)=t^4$ satisfies the second condition with
$$\alpha(x,y)=\frac{4y^3(y-x)}{(x-y)^2}>0,$$
but it is not strongly quasiconvex on any open interval containing $0$: otherwise, by \cite[Corollary 9]{hadjisavvas2026characterizations}, we would have $t^4\geq\frac{\gamma}{4}t^2$ for all $t$ in that interval, whereas $t^4<\frac{\gamma}{4}t^2$ for all sufficiently small $t\neq 0$.
\end{remark}

\section{Consequences}\label{Sec. Conseq}

In this section, we apply Theorem \ref{teo. Clarke} to the problem
\begin{equation}\label{eq_PO}
\begin{aligned}
\min_{x}\quad & h(x)\\
\text{s.t.}\quad & x\in K,
\end{aligned}
\tag{P}
\end{equation}
where $C\subseteq \mathbb{R}^n$ is an open convex set, $h:C\to \mathbb{R}$ is locally Lipschitz and strongly quasiconvex with modulus $\gamma >0$, and $K\subseteq C$ is nonempty, closed and convex. No convexity, differentiability or Clarke regularity assumption is imposed on $h$. We also study the relation between proximal points and strong subdifferentials.

\subsection{Optimality and the subgradient method}\label{subsec_optim}

We denote by $N_K(x):=\{v\in \mathbb{R}^n: \langle v, y-x \rangle\leq 0 \text{ for all }y \in K\}$ the normal cone to $K$ at $x\in K$, by $P_K$ the metric projection onto $K$, and by $\dist(y,S)$ the distance from a point $y$ to a set $S$. Since $h$ is locally Lipschitz, its restriction to $K$ is continuous, in particular lower semicontinuous, and strongly quasiconvex, so that \eqref{eq_PO} has a unique solution $\bar x$ \cite[Corollary 3]{lara2022strongly}. Since $h(\bar x)\leq h(x)$ for all $x\in K$, the implication (a) $\Rightarrow$ (c) of Theorem \ref{teo. Clarke} yields, for all $x\in K$ and all $z\in\partial^{\mathcal{C}}h(x)$,
\begin{equation}\label{Cons. eq1}
    \langle z,x-\bar x\rangle \geq \frac{\gamma}{2}\|x-\bar x\|^2.
\end{equation}

\begin{corollary}\label{Cons. coro_opt1}
A point $\bar x\in K$ is a solution of \eqref{eq_PO} if and only if $0\in \partial^\mathcal{C} h(\bar x)+N_K(\bar x)$. Moreover, for $\bar x\in C$, we have $0\in \partial^\mathcal{C} h(\bar x)$ if and only if $\bar x$ is the unique minimizer of $h$ on $C$.
\end{corollary}
\begin{proof}
($\Rightarrow$) If $\bar x$ is a minimizer of $h$ on $K$, since $h$ is Lipschitz near $\bar x$, we have $0\in \partial^\mathcal{C}h(\bar x)+N_K(\bar x)$ \cite[Corollary to Proposition 2.4.3]{clarke1990optimization}, that is, there exists $z\in\partial^{\mathcal C}h(\bar x)$ with $-z\in N_K(\bar x)$. In this inclusion, the Clarke normal cone coincides with $N_K(\bar x)$, since $K$ is convex \cite[Proposition 2.4.4]{clarke1990optimization}.

($\Leftarrow$) Let $z\in \partial^\mathcal{C} h(\bar x)$ with $-z\in N_K(\bar x)$, that is, $\langle z,y-\bar x\rangle\geq 0$ for all $y\in K$. If there existed $y\in K\setminus \{ \bar x\}$ with $h(y)\leq h(\bar x)$, item (c) of Theorem \ref{teo. Clarke} would imply $ \langle z,y-\bar x\rangle \leq -\frac{\gamma}{2}\|y-\bar x\|^2<0$, which is absurd. Hence,
$h(y)>h(\bar x)$ for all $y\in K\setminus\{\bar x\}$. Thus, $\bar x$ is the unique minimizer of $h$ on $K$.

Finally, let $\bar x\in C$. If $\bar x$ is a minimizer of $h$ on $C$, then $\bar x$ is a local minimizer, since $C$ is open, and $0\in\partial^{\mathcal{C}}h(\bar x)$ \cite[Proposition 2.3.2]{clarke1990optimization}. Conversely, if $0\in\partial^{\mathcal{C}}h(\bar x)$ and there existed $y\in C\setminus \{\bar x\}$ with $h(y)\leq h(\bar x)$, item (c) of Theorem \ref{teo. Clarke} would imply $0 \leq -\frac{\gamma}{2}\|y-\bar x\|^2<0$, which is absurd.
\end{proof}

For $K=C=\mathbb{R}^n$, the sufficient condition of Corollary \ref{Cons. coro_opt1} coincides with \cite[Lemma 1]{vial1982strong}. For $K$ closed and convex and contained in an open convex set, the characterization in Corollary \ref{Cons. coro_opt1} is \cite[Corollary 3.1]{dudov2022sufficient}, under the additional assumption that $h$ is Clarke regular,
an assumption that is not needed in Corollary \ref{Cons. coro_opt1}. Moreover, applied to the restriction of $h$ to $K$, \cite[Corollary 9]{hadjisavvas2026characterizations} yields $h(y)-h(\bar x)\geq\frac{\gamma}{4}\|y-\bar x\|^2$ for all $y\in K$, which improves the constant $\frac{\gamma}{8}$ of \cite[(3.8)]{dudov2022sufficient}; by Remark \ref{MR. obs2}, the constant $\frac{\gamma}{4}$ is optimal.

\begin{corollary}\label{Cons. coro_opt2}
If $\bar x$ is a minimizer of $h$ on $K$, then
\begin{equation}\label{eq1.coro_opt2}
    \|x-\bar x\|\leq \frac{2}{\gamma}\dist \left(0,\partial^\mathcal{C} h(x)+N_K(x)\right)\quad \text{for all } x\in K.
\end{equation}
\end{corollary}
\begin{proof}
For $x=\bar x$, the inequality is immediate. For $x\in K\setminus \{\bar x\}$, let $z\in \partial^\mathcal{C}h(x)$ and $u\in N_K(x)$. Since $\langle u,\bar x-x\rangle \leq 0$, it follows from \eqref{Cons. eq1} that $\langle z+u,x-\bar x\rangle \geq \frac{\gamma}{2}\|x-\bar x\|^2$, and, by the Cauchy--Schwarz inequality, we obtain
$$\|z+u\|\geq \frac{\gamma}{2}\|x-\bar x\|.$$
Taking the infimum over $z$ and $u$, we obtain \eqref{eq1.coro_opt2}.
\end{proof}

\begin{remark}
The constant $\frac{2}{\gamma}$ in \eqref{eq1.coro_opt2} is optimal, that is, it cannot be replaced by $\frac{c}{\gamma}$ with $c<2$. Indeed, let $b>1$, $C=(-2,2)$, $K=[-1,1]$, $h_b(t):=t^2$ for $t\geq 0$ and $h_b(t):=bt^2$ for $t< 0$, whose exact modulus on $C$ is $\gamma_b=\frac{4\sqrt{b}}{1+\sqrt{b}}$, by Remark \ref{MR. obs2}. We have $\bar x=0$ and, for $x=t\in(0,1)$, $N_K(t)=\{0\}$ and $\partial^\mathcal{C}h_b(t)=\{2t\}$. Since $c<2$ and $\frac{1+\sqrt{b}}{\sqrt{b}}\to 1$ as $b\to\infty$, for $b$ sufficiently large we obtain
$$\frac{c}{\gamma_b}\dist \left(0,\partial^\mathcal{C} h_b(t)+N_K(t)\right)= \frac{c}{2}\,\frac{1+\sqrt{b}}{\sqrt{b}}\,t\ <\ t=|t-\bar x|. $$
Thus, \eqref{eq1.coro_opt2} with $\frac{c}{\gamma}$ in place of $\frac{2}{\gamma}$ fails at $x=t$.
\end{remark}

If $K=C=\mathbb{R}^n$ and $h$ is differentiable with $L$-Lipschitz gradient $\nabla h$, then $\partial^\mathcal{C}h(x)=\{\nabla h(x)\}$ and \eqref{eq1.coro_opt2} reduces to
$$\frac{\gamma}{2}\|x-\bar x\|\leq \|\nabla h(x)\|.$$
Combining this inequality with $h(x)-h(\bar x)\leq \frac{L}{2}\|x-\bar x\|^2$, which follows from the descent lemma \cite[Lemma 1.2.3]{Nesterov2018} since $\nabla h(\bar x)=0$, we obtain the Polyak--{\L}ojasiewicz inequality
$$\frac{\gamma^2}{2L}(h(x)-h(\bar x))\leq \|\nabla h(x)\|^2$$
of \cite[Proposition 6]{lara2025characterizations}.
Thus, \eqref{eq1.coro_opt2} extends the bound underlying this inequality to Clarke subgradients and to constrained problems, with the optimal constant.

We now consider the projected subgradient method, in which $z_k$ is an arbitrary Clarke subgradient. For strongly quasiconvex functions, this scheme was studied in \cite[Algorithm 1]{choque2024subgradient} with subgradients in the strong subdifferential $\partial^K_{\beta,\gamma}h(x_k)$.

\begin{algorithm}[H]
\caption{Projected subgradient method with Clarke subgradients}\label{alg. SPC}
 \raggedright
Step 0. Choose $x_0\in K$ and a sequence $(\alpha_k)_{k\geq 0}$ of positive numbers, and set $k=0$.\\[2pt]
Step 1. Choose $z_k\in \partial^\mathcal{C}h(x_k)$. If $z_k=0$, \textsc{stop}.\\[2pt]
Step 2. Set
\begin{equation}\label{Cons. eq2}      x_{k+1}:=P_K(x_k-\alpha_kz_k),
\end{equation}
$k:=k+1$, and go to Step 1.
\end{algorithm}

The stopping criterion is justified by Corollary \ref{Cons. coro_opt1}, since $z_k=0$ implies $0\in \partial^\mathcal{C}h(x_k)+N_K(x_k)$. Moreover, if $K$ is compact, then $\ell:=\sup\{\|z\|:~z\in\partial^\mathcal{C}h(x),~x\in K\}$ is finite, since $h$ is Lipschitz on a neighborhood of $K$ and the subgradients are bounded by the Lipschitz constant \cite[Proposition 2.1.2(a)]{clarke1990optimization}.

\begin{theorem}\label{Cons. teo1}
Let $K$ be compact, $\bar x$ the minimizer of $h$ on $K$, and $(x_k)$ the sequence generated by Algorithm \ref{alg. SPC}. Then the following statements hold:
\begin{enumerate}
\item [(a)] If $z_k=0$ for some $k$, then $x_k=\bar x$.
\item [(b)] If $\alpha_k=\alpha\in (0,\frac{1}{\gamma})$ for all $k$, then
$$\|x_k-\bar x\|^2\leq (1-\alpha\gamma)^k\|x_0-\bar x\|^2+\frac{\alpha \ell^2}{\gamma} \quad \text{for all } k.$$
\item [(c)] If $\alpha_k=\frac{2}{\gamma(k+2)}$ for all $k$, then, for all $k\geq 1$,
    $$\|x_k-\bar x\|^2\leq \frac{4 \ell^2}{\gamma^2(k+1)}\quad \quad \text{and} \quad \quad h(x_k)-h(\bar x)\leq \frac{2\ell^2}{\gamma\sqrt{k+1}}.$$
\end{enumerate}
\end{theorem}
\begin{proof} Let $a_k:=\|x_k-\bar x\|^2$. Since $P_K$ is nonexpansive, $P_K\bar x=\bar x$ and $x_k\in K$, we have
$$a_{k+1}\leq\|x_k-\bar x-\alpha_kz_k\|^2=a_k-2\alpha_k\langle z_k,x_k-\bar x\rangle +\alpha_k^2\|z_k\|^2.$$
By \eqref{Cons. eq1} and $\|z_k\|\leq \ell$, it follows that
\begin{equation}\label{Cons. teo eq1}
    a_{k+1} \leq (1-\alpha_k \gamma)a_k+\alpha_k^2\ell^2.
\end{equation}
\begin{enumerate}
    \item [(a)] If $z_k=0$, \eqref{Cons. eq1} implies $\|x_k-\bar x\|=0$, that is, $x_k=\bar x$.
    \item [(b)] Since $0< 1-\alpha\gamma<1$, iterating \eqref{Cons. teo eq1} gives $$a_k\leq (1-\alpha\gamma)^ka_0+\alpha^2\ell^2\sum_{j= 0}^{k-1}(1-\alpha\gamma)^j\leq (1-\alpha\gamma)^ka_0+\frac{\alpha\ell^2}{\gamma}.$$
    \item[(c)] Since $\alpha_k=\frac{2}{\gamma(k+2)}$, \eqref{Cons. teo eq1} reads
$$a_{k+1}\leq \frac{k}{k+2}a_k+\frac{4\ell^2}{\gamma^2(k+2)^2}.$$
Multiplying by $(k+1)(k+2)$ and using $\frac{k+1}{k+2}\leq 1$, we obtain
$$(k+1)(k+2)a_{k+1}\leq k(k+1)a_k + \frac{4\ell^2}{\gamma^2},$$
and, summing these inequalities from $0$ to $k-1$, we get
$k(k+1)a_k\leq  \frac{4\ell^2}{\gamma^2}k$, that is,
$$\|x_k-\bar x\|^2\leq \frac{4\ell^2}{\gamma^2(k+1)}.$$
If $x_k\neq \bar x$, Lebourg's mean value theorem \cite[Theorem 2.3.7]{clarke1990optimization} guarantees the existence of $u \in (\bar x,x_k)$ and $\zeta\in \partial^{\mathcal{C}}h(u)$ such that $h(x_k)-h(\bar x)=\langle \zeta, x_k-\bar x \rangle$. Since $K$ is convex, $u\in K$; hence $\|\zeta \|\leq \ell$ and
$$h(x_k)-h(\bar x)\leq \ell\|x_k-\bar x\|\leq \frac{2\ell^2}{\gamma\sqrt{k+1}},$$
an estimate that is immediate if $x_k=\bar x$.
\end{enumerate}
\end{proof}

\begin{remark}\label{Obs_Cons_teo1}
\begin{enumerate}
    \item [(i)] The step size and the bound in (c) coincide, in the deterministic case, with those of \cite[Section 3.2]{lacoste2012simpler} for strongly convex functions, and Theorem \ref{Cons. teo1} shows that they persist under strong quasiconvexity. The step size in (c) depends only on $\gamma$, and the corresponding bound is independent of the starting point $x_0$. Since this bound decreases as $\gamma$ increases, it is advantageous to use the exact modulus, which can be computed by means of Theorem \ref{teo. P}, as in Example \ref{exem. Table}.

    \item [(ii)] Unlike the elements of $\partial^K_{\beta,\gamma}h(x_k)$ used in \cite{choque2024subgradient}, Clarke subgradients can be obtained through calculus rules. For instance, if $h=\max_ih_i$ with each $h_i$ of class $C^1$, then, since each $h_i$ is regular \cite[Proposition 2.3.6(a)]{clarke1990optimization}, we have $\partial^{\mathcal{C}}h(x_k)=\conv\{\nabla h_i(x_k): \ i\in I(x_k) \}$ \cite[Proposition 2.3.12]{clarke1990optimization}, where $I(x_k)$ is the set of active indices at $x_k$. In particular, the gradient of any $h_i$ that is active at $x_k$ is a Clarke subgradient of $h$ at $x_k$.

    \item [(iii)] For Lipschitz strongly quasiconvex functions, \cite[Remark 19(i)]{lara2022strongly} suggests the use of subgradient-type methods, such as that of \cite[Section~3]{kiwiel2001convergence}, whose directions are normalized and taken in normal cones to sublevel sets. Theorem \ref{Cons. teo1} shows that non-normalized Clarke subgradients, with step size $\frac{2}{\gamma(k+2)}$, yield the rate $O(\frac{1}{k})$ in (c).
\end{enumerate}
\end{remark}

With $K=C=\mathbb{R}^n$, \eqref{Cons. eq2} is the explicit Euler discretization of the Clarke subgradient flow, for which \eqref{Cons. eq1} guarantees exponential convergence.

\begin{proposition}\label{Cons. Prop1}
Let $h: \mathbb{R}^n\to \mathbb{R}$ be locally Lipschitz and strongly quasiconvex with modulus $\gamma>0$, $\bar x$ its unique minimizer, $x_0 \in \mathbb{R}^n$, $B:=\{x: \|x-\bar x\|\leq\| x_0-\bar x\|\}$, and $\ell_0:=\sup\{\|z\|: z\in\partial^{\mathcal{C}}h(x),~x\in B\}$. Let $0<T \leq \infty$ and let $x:[0,T)\to \mathbb{R}^n$ be absolutely continuous on every compact subinterval, with
\begin{equation}\label{Cons. eq3}
     \dot x(t)\in -\partial^{\mathcal{C}}h(x(t))\quad  \text{for almost every }  t\in [0,T),\quad x(0)=x_0.
\end{equation}
Then, for all $ t\in [0,T)$,
$$\|x(t)-\bar x\|\leq e^{-\frac{\gamma t}{2}}\|x_0-\bar x\|,\quad \quad h(x(t))-h(\bar x)\leq \ell_0\,e^{-\frac{\gamma t}{2}}\|x_0-\bar x\|.$$
Moreover, for every $x_0\in\mathbb{R}^n$, \eqref{Cons. eq3} admits at least one solution defined on $[0,\infty)$.
\end{proposition}

\begin{proof} The function $r(t):=\|x(t)-\bar x\|^2$ is absolutely continuous on every compact subinterval of $[0,T)$. Since $-\dot x(t)\in \partial^\mathcal{C}h(x(t))$ for almost every $t$, \eqref{Cons. eq1} with $K=\mathbb{R}^n$ yields
$$\dot r(t)=2\langle \dot x(t), x(t)-\bar x \rangle\leq -\gamma r(t) \quad \text{for almost every } t\in [0,T).$$
Hence, $\frac{d}{dt}\left( e^{\gamma t}r(t)\right)\leq 0$ for almost every $t$ and, since the function $t\mapsto e^{\gamma t}r(t)$ is absolutely continuous on compact subintervals, it is nonincreasing. Thus, $e^{\gamma t} r(t)\leq r(0)$, that is,
$$\|x(t)-\bar x\|\leq e^{-\frac{\gamma t}{2}}\|x_0-\bar x\|\quad \text{for all } t\in [0,T).$$
In particular, $x(t)\in B$ for all $t\in [0,T)$. Since $B$ is compact, $h$ is Lipschitz on a neighborhood of $B$, and the Clarke subgradients are bounded by the Lipschitz constant \cite[Proposition 2.1.2(a)]{clarke1990optimization}, so that $\ell_0<\infty$. If $x(t)\neq \bar x$, by Lebourg's mean value theorem \cite[Theorem 2.3.7]{clarke1990optimization}, there exist $u \in \{\bar x+s(x(t)-\bar x): 0< s<1\}\subseteq B$ and
$\zeta \in \partial^\mathcal{C}h(u)$ such that
$$h(x(t))-h(\bar x)=\langle \zeta, x(t)-\bar x\rangle.$$
Since $\|\zeta \|\leq \ell_0$, it follows that
$$h(x(t))-h(\bar x)\leq \ell_0\|x(t)-\bar x\|\leq \ell_0\,e^{-\frac{\gamma t}{2}} \|x_0-\bar x\|,$$
an inequality that holds trivially if $x(t)=\bar x$.

For the existence, let $B_1:=\{y: \dist(y,B)\leq 1\}$, $M:=\sup \{\|z\|: z\in \partial^\mathcal{C}h(y),~y\in B_1\}$, which is finite by the same argument, and $a:=\frac{1}{M+1}$. The map $-\partial^\mathcal{C}h$ is upper semicontinuous, with nonempty, convex and compact values \cite[Propositions 2.1.2(a) and 2.1.5(d)]{clarke1990optimization}. For $y_0\in B$, we have $\{y:\ \|y-y_0\|\leq 1 \}\subseteq B_1$, so that
 $$\|z\|\leq M \quad \text{whenever } z\in\partial^{\mathcal C}h(y) \text{ and }  \|y-y_0\|\le1.$$
Hence, there exists a solution of
 \eqref{Cons. eq3} on $[0,a]$ with $x(0)=y_0$ \cite[Chap. 2, Sec. 2.1, Theorem 3]{aubin1984differential}. By the first part, $\|x(a)-\bar x\|\leq \|y_0-\bar x\|$, that is, $x(a)\in B$. Since \eqref{Cons. eq3} is autonomous, restarting at $x(a), x(2a),\ldots$ and concatenating, we obtain a solution on $[0,\infty)$.
\end{proof}

Proposition \ref{Cons. Prop1} extends \cite[Theorem 15(b),(c)]{lara2025characterizations}, where $h$ is of class $C^1$, to locally Lipschitz functions, with the same rate $e^{-\frac{\gamma }{2}t}$. In the estimate for the values of $h$, the constant $\ell_0$ depends only on $h$ and $x_0$.
As in \cite{lara2025characterizations}, the gradient (or subgradient) is not assumed to be Lipschitz,
and the estimates apply to every solution of \eqref{Cons. eq3}, which is not necessarily unique. The monotonicity of $t\mapsto h(x(t))$, obtained in \cite[Theorem~15(a)]{lara2025characterizations} in the $C^1$ case, is discussed in Remark \ref{Obs.aberto}.

\begin{remark}\label{obs_Cons. Prop1}
For $t \geq 0$ and $x_0\in\mathbb{R}^n$, let
$$
G(t,x_0):=\{x(t):~ x\text{ is a solution of \eqref{Cons. eq3} on }[0,\infty)\text{ with }x(0)=x_0\}.
$$
By Proposition~\ref{Cons. Prop1}, $G(t,x_0)\neq\emptyset$ and, since \eqref{Cons. eq3} is autonomous, $G(0,x_0)=\{x_0\}$ and $G(t_1+t_2,x_0)\subseteq G(t_1,G(t_2,x_0))$, that is, $G$ is a multivalued semiflow in the sense of \cite[Definition 1 and Remark~1]{melnik1998attractors}. Proposition \ref{Cons. Prop1} gives $G(t,\bar x)=\{\bar x\}$ and, for every bounded $D\subseteq\mathbb{R}^n$,
\[
\sup_{y\in G(t,D)}\|y-\bar x\|\leq e^{- \frac{\gamma}{2}t}\sup_{x_0\in D}\|x_0-\bar x\|,
\]
so that $\{\bar x\}$ is a global attractor of $G$ \cite[Definition 6]{melnik1998attractors}, with exponential attraction. The Lyapunov function used is $V(x)=\|x-\bar x\|^2$, not $h$. Unlike the strongly convex case, in which the strong monotonicity of the subdifferential brings any two trajectories exponentially close to each other, \eqref{Cons. eq1} compares each trajectory only with $\bar x$, which is compatible with the possible non-uniqueness of solutions.
\end{remark}

The following example illustrates Theorem \ref{Cons. teo1} with a nonconvex function, given by the maximum of two smooth functions, whose minimizer lies in the set where it is not differentiable. By Remark \ref{Obs_Cons_teo1}(ii), a Clarke subgradient of this function is the gradient of any of the functions that is active at the point.

\begin{example}\label{Exem.alg. SPC}
Let $f:\mathbb{R} \to \mathbb{R}$ be given by $f(r)=r^2+2\sin^2r$, let $a=(0.6, 0)$, $b=(-0.6,0.3)$, $g_a(x)=f(\|x-a\|)$, $g_b(x)=f(\|x-b\|)$, and $h=\max\{g_a,g_b\}$ on $C=(-\frac{3}{2}, \frac{3}{2})\times (-\frac{3}{2}, \frac{3}{2})$, with $K=[-1,1]\times[-1,1]$.

The set $C$ is open and convex, and $K\subseteq C$ is compact and convex. Since $f(r)=r^2+1-\cos(2r)$ is analytic and even, the functions $g_a$ and $g_b$ are of class $C^1$; hence, $h$ is locally Lipschitz. As for strong quasiconvexity, since $h$ is not differentiable, we verify \eqref{eq. grad} for $g_a$ and $g_b$. Since $f^\prime(r)=2r+2\sin2r>0$ for $r>0$,
$$g_a(x)\leq g_a(y) \quad  \iff \quad \|x-a\|\leq \|y-a\|.$$
If $y=a$, this condition implies $x=a$, and \eqref{eq. grad} is immediate.
For $\rho:=\|y-a\|>0$, we have $$\nabla g_a(y)=\frac{f^\prime(\rho)}{\rho}(y-a).$$ If $\|x-a\|\leq \|y-a\|$ then, with $u:=y-a$ and $v:=x-a$,
$$\langle  u,u-v\rangle-\frac{1}{2}\|u-v\|^2=\frac{1}{2}(\|u\|^2-\|v\|^2)\geq 0.$$
Hence,
$$\langle \nabla g_a(y),y-x\rangle =\frac{f^\prime (\rho)}{\rho}\langle u,u-v\rangle \geq \frac{1}{2} \frac{f^\prime(\rho)}{\rho}\|y-x\|^2.$$
Thus, \eqref{eq. grad} is satisfied on $C$ with $$\gamma=\inf_{r>0}\frac{f^\prime(r)}{r}=2+4\cos u^*\approx 1.1311,$$ where $u^*\approx 4.4934$ is the first positive solution of $\tan u=u$. By Corollary \ref{MR cor1}, $g_a$ is strongly quasiconvex on $C$ with modulus $\gamma$; the same argument applies to $g_b$. Consequently, $h=\max\{g_a,g_b\}$ is also strongly quasiconvex with modulus $\gamma$ \cite[Proposition 11(c)]{grad2025strongly}.

The function $h$ is not convex on $K$, since, for $p=(-0.5,0)$ and $q=(-1,0)$, we get $$h\left(\frac{p+q}{2}\right)\approx 3.727>\frac{h(p)+h(q)}{2}\approx 3.678.$$ Moreover, since $\max \{\|x-a\|,\|x-b\|\}\geq \frac{1}{2}\|a-b\|$ and $f$ is increasing on $[0,\infty)$, we have $$h(x)\geq f\left(\frac{1}{2}\|a-b\|\right),$$ with equality only at $\bar x=\frac{1}{2}(a+b)=(0,0.15)$. Thus, $\bar x$ is the minimizer of $h$ on $K$ and lies on the perpendicular bisector of the segment $[a,b]$, on which $g_a=g_b$ and $\nabla g_a \neq \nabla g_b$, so that $h$ is not differentiable at $\bar x$.

As for the method, the Clarke subgradients of $h$ on $K$ are convex combinations of the gradients of the active functions \cite[Proposition~2.3.12]{clarke1990optimization}, whose norms are of the form $f^\prime(\rho)$ with $\rho\leq 2.07$. Hence,
$$\ell\leq \max_{0 \leq \rho\leq 2.07}f^\prime(\rho)
=f^\prime\left(\frac{\pi}{3}\right)=\frac{2\pi}{3}+\sqrt{3}\approx 3.83.$$
Since \eqref{Cons. eq1} remains true with any $\gamma^\prime \in(0,\gamma]$ in place of $\gamma$, Theorem \ref{Cons. teo1}(c) applies with $\gamma^\prime=1.1$. Taking, at each iteration, $\alpha_k:=\frac2{\gamma^\prime(k+2)}$ and
$$z_k:=\begin{cases}\nabla g_a(x_k), & \text{if } g_a(x_k)\ge g_b(x_k),\\ \nabla g_b(x_k), & \text{if } g_a(x_k)<g_b(x_k),\end{cases}$$
that is, the gradient of a function that is active at $x_k$, we have $z_k\in\partial^\mathcal{C}h(x_k)$, and every sequence generated by Algorithm \ref{alg. SPC} satisfies
 $$\|x_k-\bar x\|^2\leq \frac{4\ell^2}{{\gamma^\prime}^2(k+1)}<\frac{49}{k+1}.$$
Figure \ref{fig:ex4.9-metodo} shows the first iterations starting from $x_0=(1,-1)$ and from $x^\prime_0=(-1,1)$, as well as the error $\|x_k-\bar x\|^2$, which remains below this bound.
\end{example}

\definecolor{corA}{RGB}{0,90,170}
\definecolor{corB}{RGB}{200,60,30}
\definecolor{corN}{RGB}{40,40,40}
\definecolor{corL}{RGB}{20,120,80}

% \lente{r}{style}: fills the lens B[a,r] \cap B[b,r]
\newcommand{\lente}[2]{%
  \begin{scope}
    \clip (0.6,0) circle (#1);
    \fill[#2] (-0.6,0.3) circle (#1);
  \end{scope}}
% \bordalente{r}{style}: draws the level curve h = f(r)
\newcommand{\bordalente}[2]{%
  \begin{scope}
    \clip (-0.6,0.3) circle (#1);
    \draw[#2] (0.6,0) circle (#1);
  \end{scope}
  \begin{scope}
    \clip (0.6,0) circle (#1);
    \draw[#2] (-0.6,0.3) circle (#1);
  \end{scope}}

% Common background: sublevel sets of h in C, the square K, and the points a and b
\newcommand{\fundosubnivel}{%
  \begin{scope}
    \clip (-1.5,-1.5) rectangle (1.5,1.5);
    \fill[corL!4] (-1.5,-1.5) rectangle (1.5,1.5);
    \foreach \r/\c in {2.4/8,2.1/13,1.8/18,1.5/24,1.2/31,0.95/39,0.75/48}{%
      \lente{\r}{corL!\c}%
      \bordalente{\r}{corL!70!black, line width=0.35pt}%
    }
  \end{scope}
  \draw[line width=0.6pt] (-1.5,-1.5) rectangle (1.5,1.5);
  \draw[line width=0.9pt, dashed] (-1,-1) rectangle (1,1);
  \node[below right, font=\footnotesize] at (-1.5,1.5) {$C$};
  \node[above left, font=\footnotesize] at (1,-1) {$K$};
  \fill[corA] (0.6,0) circle (1.2pt) node[below right=-1pt, font=\footnotesize] {$a$};
  \fill[corB] (-0.6,0.3) circle (1.2pt) node[above left=-1pt, font=\footnotesize] {$b$};
}

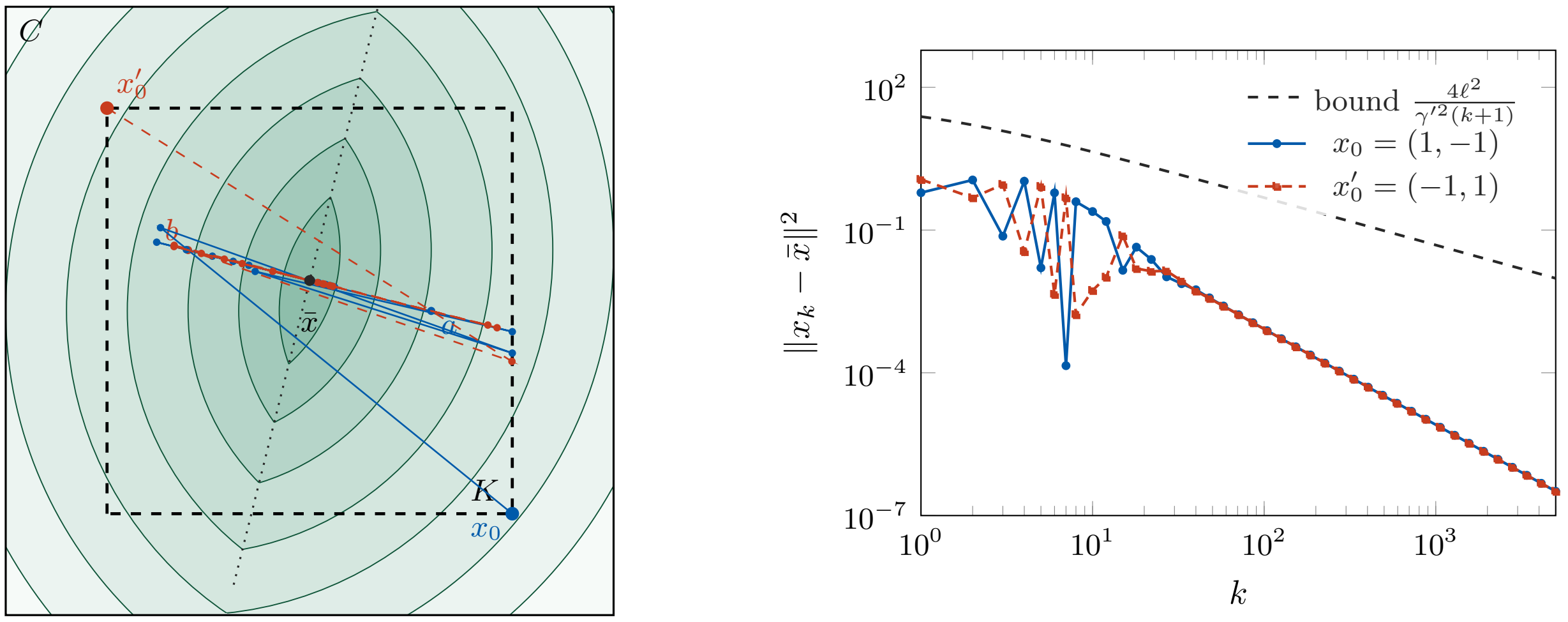
\begin{figure}[H]
\centering
\begin{minipage}[b]{0.48\textwidth}\centering
\begin{tikzpicture}[x=2.1cm,y=2.1cm]
  \fundosubnivel
  \begin{scope}
    \clip (-1.5,-1.5) rectangle (1.5,1.5);
    \draw[corN, line width=0.6pt, dotted]
      (-0.375,-1.35) -- (0.4125,1.8);
  \end{scope}
  % x0 = (1,-1): first 40 iterations
  \draw[corA, line width=0.5pt] plot coordinates {
    (1.0000,-1.0000) (-0.7360,0.4105) (1.0000,-0.2079) (-0.2690,0.1949)
    (1.0000,-0.1025) (0.1228,0.1182) (-0.7548,0.3389) (0.0116,0.1472)
    (-0.6088,0.3022) (0.0441,0.1390) (-0.4801,0.2700) (0.0802,0.1300)
    (-0.3776,0.2444) (0.1038,0.1240) (-0.3003,0.2251) (0.1159,0.1210)
    (-0.2438,0.2110) (0.1198,0.1200) (-0.2029,0.2007) (0.1188,0.1203)
    (-0.1730,0.1932) (0.1151,0.1212) (-0.1506,0.1877) (0.1101,0.1225)
    (-0.1335,0.1834) (0.1045,0.1239) (-0.1201,0.1800) (0.0989,0.1253)
    (-0.1093,0.1773) (0.0936,0.1266) (-0.1005,0.1751) (0.0885,0.1279)
    (-0.0931,0.1733) (0.0838,0.1291) (-0.0868,0.1717) (0.0794,0.1301)
    (-0.0813,0.1703) (0.0754,0.1311) (-0.0766,0.1691) (0.0718,0.1321)
    (-0.0724,0.1681)};
  \foreach \p in {(1.0000,-1.0000),(-0.7360,0.4105),(1.0000,-0.2079),
    (-0.2690,0.1949),(1.0000,-0.1025),(0.1228,0.1182),(-0.7548,0.3389),
    (0.0116,0.1472),(-0.6088,0.3022),(0.0441,0.1390),(-0.4801,0.2700),
    (0.0802,0.1300),(-0.3776,0.2444),(0.1038,0.1240),(-0.3003,0.2251)}
    \fill[corA] \p circle (1.1pt);
  % x0' = (-1,1): first 40 iterations
  \draw[corB, line width=0.5pt, dashed] plot coordinates {
    (-1.0000,1.0000) (0.9982,-0.2489) (-0.6686,0.3235) (0.9249,-0.0829)
    (-0.1822,0.1951) (0.8794,-0.0697) (0.0647,0.1339) (-0.6696,0.3174)
    (0.0396,0.1401) (-0.5346,0.2836) (0.0711,0.1322) (-0.4211,0.2553)
    (0.0981,0.1255) (-0.3330,0.2333) (0.1139,0.1215) (-0.2678,0.2169)
    (0.1204,0.1199) (-0.2204,0.2051) (0.1210,0.1198) (-0.1859,0.1965)
    (0.1180,0.1205) (-0.1604,0.1901) (0.1132,0.1217) (-0.1411,0.1853)
    (0.1077,0.1231) (-0.1261,0.1815) (0.1020,0.1245) (-0.1142,0.1786)
    (0.0964,0.1259) (-0.1045,0.1761) (0.0911,0.1272) (-0.0965,0.1741)
    (0.0862,0.1285) (-0.0897,0.1724) (0.0816,0.1296) (-0.0839,0.1710)
    (0.0775,0.1306) (-0.0788,0.1697) (0.0736,0.1316) (-0.0744,0.1686)
    (0.0701,0.1325)};
  \foreach \p in {(-1.0000,1.0000),(0.9982,-0.2489),(-0.6686,0.3235),
    (0.9249,-0.0829),(-0.1822,0.1951),(0.8794,-0.0697),(0.0647,0.1339),
    (-0.6696,0.3174),(0.0396,0.1401),(-0.5346,0.2836),(0.0711,0.1322),
    (-0.4211,0.2553),(0.0981,0.1255),(-0.3330,0.2333),(0.1139,0.1215)}
    \fill[corB] \p circle (1.1pt);
  \fill[corA] (1,-1) circle (2pt) node[below left=-1pt, font=\footnotesize] {$x_0$};
  \fill[corB] (-1,1) circle (2pt) node[above right=-1pt, font=\footnotesize] {$x_0'$};
  \fill[corN] (0,0.15) circle (1.6pt);
  \node[font=\footnotesize, anchor=north] at (0,0.06) {$\bar x$};
\end{tikzpicture}
\\[2pt]{\small (a)}
\end{minipage}\hfill
\begin{minipage}[b]{0.48\textwidth}\centering
\begin{tikzpicture}
\begin{loglogaxis}[
  width=\linewidth, height=6.4cm,
  xmin=1, xmax=5000, ymin=1e-7, ymax=600,
  xlabel={$k$}, ylabel={$\|x_k-\bar x\|^2$},
  tick label style={font=\footnotesize}, label style={font=\small},
  legend pos=north east,
  legend style={font=\footnotesize, draw=none, fill=white, fill opacity=0.85}]
  \addplot[domain=1:5000, samples=60, thick, corN, dashed] {48.4023/(x+1)};
  \addlegendentry{bound $\frac{4\ell^2}{\gamma'^2(k+1)}$}
  \addplot[corA, thick, mark=*, mark size=0.9pt] coordinates {
    (1,6.0951e-01) (2,1.1281e+00) (3,7.4392e-02) (4,1.0638e+00)
    (5,1.6105e-02) (6,6.0547e-01) (7,1.4119e-04) (8,3.9385e-01)
    (10,2.4492e-01) (12,1.5150e-01) (15,1.4268e-02) (18,4.3742e-02)
    (22,2.4102e-02) (27,1.0400e-02) (33,7.4575e-03) (40,5.5634e-03)
    (48,3.7635e-03) (58,2.5344e-03) (71,1.6624e-03) (86,1.1302e-03)
    (104,7.6816e-04) (126,5.2078e-04) (153,3.5173e-04) (186,2.3723e-04)
    (226,1.6025e-04) (274,1.0878e-04) (333,7.3506e-05) (404,4.9863e-05)
    (490,3.3852e-05) (595,2.2934e-05) (722,1.5562e-05) (876,1.0564e-05)
    (1063,7.1696e-06) (1290,4.8659e-06) (1565,3.3047e-06) (1899,2.2437e-06)
    (2305,1.5225e-06) (2797,1.0337e-06) (3395,7.0152e-07) (4120,4.7627e-07)
    (5000,3.2334e-07)};
  \addlegendentry{$x_0=(1,-1)$}
  \addplot[corB, thick, densely dashed, mark=square*, mark size=0.9pt] coordinates {
    (1,1.1556e+00) (2,4.7709e-01) (3,9.0958e-01) (4,3.5220e-02)
    (5,8.2163e-01) (6,4.4442e-03) (7,4.7637e-01) (8,1.6657e-03)
    (10,5.3719e-03) (12,1.0219e-02) (15,7.6194e-02) (18,1.5545e-02)
    (22,1.3618e-02) (27,1.3864e-02) (33,8.5495e-03) (40,5.2205e-03)
    (48,3.6533e-03) (58,2.5015e-03) (71,1.6713e-03) (86,1.1277e-03)
    (104,7.6746e-04) (126,5.2059e-04) (153,3.5178e-04) (186,2.3722e-04)
    (226,1.6024e-04) (274,1.0877e-04) (333,7.3506e-05) (404,4.9863e-05)
    (490,3.3852e-05) (595,2.2934e-05) (722,1.5562e-05) (876,1.0564e-05)
    (1063,7.1696e-06) (1290,4.8659e-06) (1565,3.3047e-06) (1899,2.2437e-06)
    (2305,1.5225e-06) (2797,1.0337e-06) (3395,7.0152e-07) (4120,4.7627e-07)
    (5000,3.2334e-07)};
  \addlegendentry{$x_0'=(-1,1)$}
\end{loglogaxis}
\end{tikzpicture}
\\[2pt]{\small (b)}
\end{minipage}
 \caption{Algorithm \ref{alg. SPC} applied to Example \ref{Exem.alg. SPC}. (a) Sublevel sets of $h$, the perpendicular bisector of $[a,b]$ (dotted), and the first 40 iterations starting from $x_0=(1,-1)$ (solid line) and from $x_0'=(-1,1)$ (dashed line). (b) Error $\|x_k-\bar x\|^2$ in log-log scale, and the bound $\frac{4\ell^2}{\gamma'^2(k+1)}\approx\frac{48.4}{k+1}$ of Theorem \ref{Cons. teo1}(c).}
\label{fig:ex4.9-metodo}
\end{figure}

\subsection{Strong subdifferentials and proximal points}\label{subsec_SFPP}

Let $K\subseteq \mathbb{R}^n$ be a closed convex set, $h:K\to \mathbb{R}$, $\beta>0$ and $\gamma\geq0$. Kabgani and Lara \cite[Definition~5]{kabgani2022strong} define the strong subdifferential $\partial_{\beta,\gamma}^Kh(x)$ as the set of all $\xi\in \mathbb{R}^n$ such that
\begin{equation}\label{Cons. eq4}
    \max\{h(y),h(x)\}\geq h(x)+\frac{\lambda}{\beta}\langle \xi,y-x\rangle + \frac{\lambda}{2}\left( \gamma-\frac{\lambda}{\beta}-\lambda\gamma\right)\|y-x\|^2
\end{equation}
for all $y\in K$ and $\lambda\in [0,1]$.
With $\Prox_{\beta h}( K,z):=\operatorname{argmin}_{y\in K}\{h(y)+\frac{1}{2\beta}\|y-z\|^2\}$, if $h$ is lower semicontinuous and strongly quasiconvex with modulus $\gamma$, then, for all $x,z\in K$,
$$x\in \Prox_{\beta h}( K,z)\quad \Longrightarrow\quad z-x\in \partial_{\beta,\gamma}^Kh(x) $$
\cite[Proposition 36]{kabgani2022strong}, and the converse is pointed out as an open problem in \cite[Remark 38]{grad2025strongly}. Proposition \ref{Cons. Prop2} shows that the converse does not hold in general, even for strongly convex functions, but that it does hold when $h$ is replaced by its truncation $h_x:=\max\{h,h(x)\}$.

\begin{proposition}\label{Cons. Prop2}
    Let $x\in K$, $z\in \mathbb{R}^n$, and let $h_x:K\to \mathbb{R}$ be given by $h_x(y):=\max\{h(y),h(x)\}$. Then the following statements hold:
    \begin{enumerate}
        \item [(a)] If $z-x\in \partial_{\beta,\gamma}^Kh(x)$, then $x\in \Prox_{\beta h_x}(K,z)$.
        \vspace{0.2cm}
        \item [(b)] The converse of \cite[Proposition 36]{kabgani2022strong} does not hold in general, even for strongly convex functions. For $K=\mathbb{R}$, $h(y)=y^2$, $\beta=\gamma=1$, $x=1$ and $z=\frac{5}{2}$, we have $z-x\in \partial_{1,1}^{\mathbb{R}}h(1)$, but $\Prox_h\left(\mathbb{R}, \frac{5}{2}\right)=\left\{\frac{5}{6}\right\}$.
    \end{enumerate}
\end{proposition}
\begin{proof}(a) Taking $\xi=z-x$ and $\lambda=1$ in \eqref{Cons. eq4}, we obtain
\begin{equation*}
    \max\{h(y),h(x)\}\geq h(x)+\frac{1}{\beta}\left( \langle z-x,y-x\rangle -\frac{1}{2} \|y-x\|^2 \right) \quad \text{for all }y\in K.
\end{equation*}
Since $\|z-y\|^2=\|z-x\|^2-2\langle z-x,y-x\rangle +\|y-x\|^2$, $\max\{h(y),h(x)\}=h_x(y)$ and $h(x)=h_x(x)$, the previous inequality is equivalent to
 $$ h_x(y)+\frac{1}{2\beta}\|z-y\|^2\geq h_x(x)+\frac{1}{2\beta} \|z-x\|^2 \quad \text{for all } y\in K.$$
Therefore, $x\in \Prox_{\beta h_x}(K,z)$.

(b) The set $K=\mathbb{R}$ is closed and convex, and $h(y)=y^2$ is continuous and strongly convex with modulus $2$, since
$$(u+t(v-u))^2=(1-t)u^2+tv^2-t(1-t)(v-u)^2,\quad  u,v\in \mathbb{R}, ~ t\in [0,1].$$
Hence, $h$ is strongly quasiconvex with modulus $2$, and therefore also with modulus $\gamma=1$, so that the assumptions of \cite[Proposition 36]{kabgani2022strong} are satisfied.
Setting $d:=y-1$ and $\xi=z-x=\frac{3}{2}$, condition \eqref{Cons. eq4} reads
$$ \max\{ d^2+2d,0\} \geq \frac{3}{2}\lambda d+\frac{\lambda}{2} (1-2\lambda)d^2\quad \text{for all } d\in \mathbb{R},~\lambda\in [0,1].$$
If $d\geq 0$, the difference between the two sides is
$$d\left(2-\frac{3}{2}\lambda\right)+d^2\left( \left(\lambda-\frac{1}{4}\right)^2+\frac{15}{16}\right)\geq 0.$$
If $-2\leq d<0$, the left-hand side vanishes and the right-hand side is
$$\lambda d\left( \frac{3}{2}+\frac{1}{2}(1-2\lambda)d\right)\leq0, $$
since $\lambda d\leq 0$ and $(1-2\lambda)d\geq-2$.
If $d<-2$, dividing by $d$, the inequality is equivalent to
$$-d\left(1-\frac{\lambda}{2}+\lambda^2\right)\geq 2-\frac{3}{2}\lambda,$$
which is true, since $-d>2$ and $2\left(1-\frac{\lambda}{2}+\lambda^2\right)= 2-\lambda+2\lambda^2 \geq 2-\frac{3}{2}\lambda$.

Thus, $z-x\in \partial_{1,1}^\mathbb{R}h(1)$. On the other hand, $ y\mapsto y^2+\frac{1}{2}(y-\frac{5}{2})^2$ is strictly convex, with derivative $3y-\frac{5}{2}$, which vanishes at $y=\frac{5}{6}$. Hence,
$$\Prox_h\left(\mathbb{R},\frac{5}{2}\right)=\left\{\frac{5}{6}\right\}\quad \text{and} \quad x=1\notin \Prox_h\left(\mathbb{R},\frac{5}{2}\right).$$
\end{proof}

\begin{remark}\label{Obs.aberto}
Two questions related to the results of this section remain open.
\begin{enumerate}
\item[(i)] Proposition \ref{Cons. Prop1} does not assert that $t\mapsto h(x(t))$ is nonincreasing. This monotonicity holds when $h$ is Clarke regular, but it may fail for general Lipschitz functions, even in $\mathbb{R}^2$ \cite[Section 1 and Proposition 3]{daniilidis2020pathological}. Does it hold for every locally Lipschitz strongly quasiconvex function?

\item [(ii)] By Proposition \ref{Cons. Prop2}, the inclusion $z-x\in\partial^K_{\beta,\gamma}h(x)$ implies $x\in\Prox_{\beta h_x}(K,z)$, but it does not imply, in general, $x\in\Prox_{\beta h}(K,z)$. For which classes of strongly quasiconvex functions does this second implication hold?
\end{enumerate}
\end{remark}

\section{Conclusions}\label{Sec. Concl}
We have shown that, for functions that are continuous along segments of a convex set, not necessarily open, the first-order inequality for the lower Dini derivative, or even for the upper one, characterizes strong quasiconvexity with the same modulus. For locally Lipschitz functions, the same characterization holds with Clarke subgradients in Banach spaces, and both characterizations rest on an elementary lemma in one variable. Since the modulus determines the constants in optimization estimates, these results lead to an optimality condition without regularity assumptions, an error bound with optimal constant, and convergence rates with the exact modulus for the projected subgradient method and for the associated flow.

\section*{Declarations}

\textbf{Funding.} J.C. Huaccha Neyra was partially supported by the São Paulo Research Foundation (FAPESP), grant \#2025/25292-1.

\textbf{Competing interests.} The authors have no competing interests to declare that are relevant to the content of this article.

\textbf{Ethics approval and consent to participate.} Not applicable.

\textbf{Consent for publication.} Not applicable.

\textbf{Data availability.} Data sharing is not applicable to this article, as no datasets were generated or analysed during the current study.

\textbf{Author contributions.} All authors contributed equally to the conception of the work, the proofs of the results and the writing of the manuscript. All authors read and approved the final manuscript.

\textbf{Use of AI tools.} During the preparation of this work, the authors used Claude (Anthropic) to improve the language, readability and clarity of the text, as well as to assist in writing the TikZ code for Figure \ref{fig:ex4.9-metodo}. After using this tool, the authors reviewed and edited the content as needed and take full responsibility for the content of this article.

%%-----------------------------
%%      bibliography
%%-----------------------------
\bibliographystyle{abbrv}
\bibliography{Bibliografia}

\end{document}